\documentclass[graybox]{svmult}

\usepackage{mathptmx}       %
\usepackage{helvet}         %
\usepackage{courier}        %
\usepackage{type1cm}        %
\usepackage{makeidx}         %
\usepackage{graphicx}        %
\usepackage{multicol}        %
\usepackage[bottom]{footmisc}%

\usepackage{amsmath}
\usepackage{amsfonts}
\usepackage{enumitem}
\usepackage{makeidx}
\usepackage{url}

\usepackage{amssymb}
\renewcommand{\qed}{\hfill$\blacksquare$}

\usepackage{mathtools}
\mathtoolsset{showonlyrefs}

\makeindex             %

\newcommand{\abs}[1]{\left\lvert#1\right\rvert}

\newcommand{\barpow}[1]{\left\lfloor#1\right\rceil}
\newcommand{\sign}[1]{\mbox{sign}(#1)}

\newcommand{\real}[1]{\mbox{Re}(#1)}

\begin{document}
\title*{Designing controllers with predefined convergence-time bound using bounded time-varying gains}
\author{Rodrigo Aldana-López, Richard Seeber, Hernan Haimovich and David~Gómez-Gutiérrez}
\authorrunning{R. Aldana-López et al.}
\titlerunning{Designing controllers with predefined convergence-time bound}
\institute{
\textcolor{red}{This is a preprint of the following chapter: Rodrigo Aldana-López, Richard Seeber, Hernan Haimovich and David~Gómez-Gutiérrez, ``Designing controllers with predefined convergence-time bound using bounded time-varying gains", publised in Sliding-Mode Control and Variable-Structure Systems, edited by Tiago Roux Oliveira, Leonid Fridman and Liu Hsu, 2023. Springer Nature Switzerland AG 2023, reproduced with permission of Springer Nature Switzerland AG. The final authenticated version is available online at: \url{https://doi.org/10.1007/978-3-031-37089-2_3}.}\\
Corresponding author: D. Gómez-Gutiérrez \\
Rodrigo Aldana-López \at Universidad de Zaragoza, Departamento de Informatica e Ingenieria de Sistemas (DIIS), María de Luna, s/n, 50018, Zaragoza, Spain. \email{rodrigo.aldana.lopez@gmail.com}
\and Richard Seeber \at Graz University of Technology, Christian Doppler Laboratory for Model Based Control of Complex Test Bed Systems, Institute of Automation and Control, Graz, Austria. \email{richard.seeber@tugraz.at}
\and Hernan Haimovich \at Centro Internacional Franco-Argentino de 
  Ciencias de la Informaci\'on y de Sistemas (CIFASIS)
  CONICET-UNR, 2000 Rosario, Argentina. \email{haimovich@cifasis-conicet.gov.ar}
\and David~Gómez-Gutiérrez\at Intel Corporation, Intel Labs, Av. del Bosque 1001, 45019, Zapopan, Jalisco, Mexico. \\
\at Tecnológico Nacional de México, Instituto Tecnológico José Mario Molina Pasquel y Henríquez, Unidad Académica Zapopan. Cam. Arenero 1101, 45019, Zapopan, Jalisco, Mexico.
\email{david.gomez.g@ieee.org}
}
\maketitle

\abstract{
Recently, there has been a great deal of attention in a class of controllers based on time-varying gains, called prescribed-time controllers,  that steer the system's state to the origin in the desired time, a priori set by the user, regardless of the initial condition. Furthermore, such a class of controllers has been shown to maintain a prescribed-time convergence in the presence of disturbances even if the disturbance bound is unknown.  However, such properties require a time-varying gain that becomes singular at the terminal time, which limits its application to scenarios under quantization or measurement noise. This chapter presents a methodology to design a broader class of controllers, called predefined-time controllers, with a prescribed convergence-time bound. Our approach allows designing robust predefined-time controllers based on time-varying gains while maintaining uniformly bounded time-varying gains. We analyze the condition for uniform Lyapunov stability under the proposed time-varying controllers.
}

\section{Introduction}
\label{sec:Intro}

Stabilizing a system in finite time in the presence of disturbance is one of the main features of sliding mode control~\cite{Shtessel2014ObservationObservers}. However, in finite-time stability, the convergence time may be an unbounded function of the system's initial condition. Thus, knowledge on the region of admissible initial conditions may be needed to deal with scenarios under time constraints. Time constraints are often present, for instance, in fault detection, isolation, and recovery schemes~\cite{Tabatabaeipour2014CalculationAnalysis}, where failing to recover from the fault on time may lead to an unrecoverable mode, or in missile guidance applications~\cite{Zarchan2012}, where the control guidance laws require stabilization in the desired time~\cite{Song2017Time-varyingTime}. 

A class of finite-time stabilization exists, called fixed-time stabilization, with a convergence-time bound independent of the initial conditions, which make it attractive to deal with time constraint. Multiple methods have been developed to obtain fixed-time stabilization, such as Lyapunov differential inequalities~\cite{Sanchez-Torres2018} and homogeneity theory~\cite{Andrieu2008}. However, not every technique allows setting a priori the desired upper bound for the convergence time, as a convergence time-bound estimate may be unknown~\cite{Andrieu2009HomogeneityDesign}. Thus, developing methods for fixed-time stabilization with a convergence-time bound defined a priori by the user has recently received a great deal of attention~\cite{Sanchez-Torres2020ASystems,Jimenez2019,Aldana-Lopez2018,Song2017Time-varyingTime}. 

On the one hand, autonomous fixed-time controllers have been explored in~\cite{Aldana-Lopez2018,Jimenez2019,Sanchez-Torres2018,Cruz-Zavala2021High-orderBi-limit,Zimenko2018}, with emphasis on estimating an upper bound for the settling time (\textit{UBST}) of the closed-loop system. Although methodologies for obtaining the least \textit{UBST} have been proposed in the literature, see e.g.~\cite{Aldana-Lopez2018,aldana2019design}, this approach has proven challenging for higher-order systems,
resulting in very conservative estimations of an \textit{UBST}~\cite{Zimenko2018}, and yielding over-engineered controllers with an unnecessarily large control magnitude.

On the other hand, prescribed-time controllers based on time-varying gains have been proposed in~\cite{Song2017Time-varyingTime,Song2019Time-varyingTime}, which have the remarkable property that, for any nonzero initial condition, its convergence time is precisely the desired one, and that no information on the disturbance bound is needed to steer the system's state to the origin. Unfortunately, the methodology requires time-varying gains that tend to infinity at the terminal time, which is problematic under quantization or measurement noise. Therefore, controllers with a predefined convergence time, taking advantage of time-varying gains while maintaining them bounded, are of great interest~\cite{Gomez2020RNC,Chitour2020StabilizationTime}. Furthermore, it is essential to analyze the uniform (with respect to time) stability property when using controllers based on time-varying gains, as the absence of uniform stability may lead to an inherent lack of robustness. However, to the best of our knowledge, such analysis is missing in the existing prescribed-time control literature.

In this chapter, we present a methodology for designing robust controllers such that the origin of its closed-loop system is fixed-time stable with a desired \textit{UBST}, i.e., predefined-time controllers. Our analysis is based on relating the closed-loop system with an auxiliary system through a time-varying coordinate change and a time-scale transformation. The methodology is motivated by an analysis of the first-order case. It is shown that applying it to a linear controller leads to a minimum energy solution. It generally allows to reduce the required control energy also when redesigning other controllers. Based on the auxiliary system's stability properties, interesting features are obtained in the closed-loop system under the proposed controller. Such an approach allows deriving a controller with the desired convergence time regardless of the initial condition, as well as predefined-time controllers with uniformly bounded time-varying gains. Finally, this methodology is complemented by studying the uniform Lyapunov stability property, providing necessary and sufficient conditions such that our methodology yields a uniformly Lyapunov stable  closed-loop system's equilibrium. 

Additionally, we show that our approach yields existing autonomous controllers as an extreme case, while the use of time-varying gains provides extra degrees of freedom for reducing the control effort.

The chapter is organized as follows: in Section~\ref{sec:motivation} we present the example of a straightforward first-order system exhibiting interesting convergence properties and from which our general strategy using time-scale transformations arises.  In Section~\ref{sec:Prelim} we provide some preliminaries on fixed-time stability and our problem of interest regarding the design of controllers with predefined convergence-time bound. In Section~\ref{sec:first_order} we provide a methodology to solve this problem in some particular cases, including first and second-order systems. We discuss some disadvantages of some prescribed-time convergence algorithms proposed in the literature, where the time-varying gains are unbounded. Finally, in Section~\ref{sec:Main} we introduce the main result of this chapter, which is the design methodology for arbitrary-order controllers with predefined convergence time-bound. In addition, we discuss the need to consider bounded time-varying gains by examining the uniform Lyapunov stability property.

\textbf{Notation:} We use boldface lower case letter for vector and boldface capital letters for matrices. The notation $\mathbf{J}:=[a_{ij}]\in\mathbb{R}^{n \times n}$ denotes a single Jordan block with zero eigenvalue, i.e., a square matrix with $a_{ij}=1$ if $j=i+1$ and $a_{ij}=0$ otherwise. The vector $\mathbf{b}_i\in\mathbb{R}^n$ denotes a vector with one in the $i$-th entry and zeros otherwise. Let $\mathbb{R}_+=\{x\in\mathbb{R}\,:\,x\geq0\}$ and $\Bar{\mathbb{R}}_+=\mathbb{R}_+\cup\{\infty\}$. 
For a function $\phi:\mathcal{I}\to\mathcal{J}$, its reciprocal $\phi(\tau)^{-1}$, $\tau\in\mathcal{I}$,  is such that $\phi(\tau)^{-1}\phi(\tau)=1$ and its inverse function $\phi^{-1}(t)$, $t\in\mathcal{J}$, is such that $\phi(\phi^{-1}(t))=t$. Given a matrix $A\in\mathbb{R}^{n\times m}$, $A^T$ represents its transpose. For a signal $y:\mathbb{R}_+\to\mathbb{R}$, $y^{(i)}(t)$ represents its $i-$th derivative with respect to time at time $t$. To denote a first-order derivative of $y(t)$, we simple use the notation $\dot{y}(t)$.

\textbf{Simulations:} Throughout the chapter, simulations are performed on OpenModelica using the Euler integration method with step size \texttt{1e-5} and tolerance \texttt{1e-6}.

\section{Motivating Example}
\label{sec:motivation}
Consider a first-order integrator 
\begin{equation}
    \label{Eq:FOI}
    \dot{x}=u+d(t)
\end{equation}
where $|d(t)|\le\Delta$ with a non-negative constant $\Delta$.
The aim is to design a feedback control law such that the origin is reached in a desired prespecified time $T_c$. 

Let us first consider the unperturbed case, i.e., where $\Delta=0$. To derive a controller, start from an auxiliary system 
\begin{equation}
\label{Eq:TauSysMotiv}
\frac{\mathrm{d}x}{\mathrm{d}\tau}=-x,
\end{equation}
written in an artificial time variable $\tau$, whose solution is
$$
x(\tau)=x_0\exp(-\tau).
$$
Our approach is to use a time-scale transformation $\tau=\varphi(t)$ such that system~\eqref{Eq:TauSysMotiv}, written with respect to the time variable $t$, reaches the origin at $t=T_c$. For this transformation to be a suitable time-scale transformation, it must be: strictly increasing, differentiable, satisfy $\lim_{t\to T_c^-}\varphi(t)=\infty$ and $\varphi(0)=0$ (a characterization of such time-scale transformations is given in~\cite{aldana2019design}). A simple example of a time-scale transformation with the above requirements is 
\begin{equation}
    \label{Eq:TimeScaleTransfMotiv}
    \tau = \varphi(t)=-\ln(1- T_c^{-1} t)
\end{equation}
whose inverse is given by 
$$t=\varphi^{-1}(\tau)=T_c\left(1-\exp(- \tau)\right).$$ Thus, the dynamics of~\eqref{Eq:TauSysMotiv} in $t$-time can be written, according to the chain rule, as
\begin{align}
\label{eq:x_solution}
\frac{\mathrm{d}x}{\mathrm{d}t}&=\left[\frac{\mathrm{d}x}{\mathrm{d}\tau}\right]_{\tau=\varphi(t)}\cdot\frac{\mathrm{d}\tau}{\mathrm{d}t}\\
&=-\frac{1}{(T_c-t)}x
\end{align}
with a solution
$$x(\varphi(t))=x_0\cdot(1- T_c^{-1} t).$$ Clearly, $$\lim_{t\to T_c^-}x(\varphi(t))=0.$$
Therefore, a controller 
\begin{equation}
\label{Eq:ControlMotiv}
u=-\frac{1}{(T_c-t)}x,
\end{equation} 
steers the state of the unperturbed integrator to the origin at a time $T_c$. 

Let us now consider the case where $\Delta\neq 0$, under the controller~\eqref{Eq:ControlMotiv}. To analyze its convergence, let us now rewrite the closed-loop system dynamics in $\tau$-time
\begin{align}
\frac{\mathrm{d}x}{\mathrm{d}\tau}&=\left[\frac{\mathrm{d}x}{\mathrm{d}t}\right]_{t=\varphi^{-1}(\tau)}\cdot\frac{\mathrm{d}t}{\mathrm{d}\tau}\\
&=\left[-\frac{1}{(T_c-t)}x+d(t)\right]_{t=\varphi^{-1}(\tau)}\cdot T_c\exp(- \tau)\\
&=-x+T_c\exp(- \tau)d(\varphi^{-1}(\tau)).
\end{align}
The solution thus satisfies
\begin{align}
    |x(\tau)|
    &=\left|x_0\exp(-\tau)+\int_{0}^{\tau}\exp(-(\tau-\xi)) T_c \exp(-\xi) d(\varphi^{-1}(\xi))\mathrm{d}\xi\right| \\
    &=\left|x_0\exp(-\tau)+T_c\exp(-\tau)\int_{0}^{\tau}d(\varphi^{-1}(\xi))\mathrm{d}\xi\right|\\
    &\leq \exp(-\tau)(|x_0|+\Delta T_c\tau).
\end{align}
Therefore, $\lim_{\tau\to\infty}x(\tau)=0$. Hence, in $t$-time, $\lim_{t\to T_c^-}x(t)=0$. 

To maintain the state at the origin after $T_c$, regardless of the disturbance, we can combine the controller~\eqref{Eq:ControlMotiv} with a sliding mode controller as follows:

\begin{equation}
\label{Eq:ControlMotiv2}
u=\left\lbrace
\begin{array}{lll}
    -\frac{1}{(T_c-t)}x & \text{for} & t\in[0,T_c) \\
     -\Delta \sign{x} & & \text{otherwise.}
\end{array}
\right.
\end{equation} 

Thus, we can summarize the following remarkable properties of this approach:
\begin{itemize}
    \item For every nonzero initial condition, the  origin is reached precisely at $T_c$, regardless of the initial conditions and without knowledge of the disturbance bound (although notice that to maintain the state at the origin after $T_c$ knowledge on the disturbance bound is required).
\end{itemize}

Unfortunately, the approach also presents the following drawback:
\begin{itemize}
    \item The time-varying gain of the controller, namely the factor $\dfrac{1}{T_c - t}$, tends to infinity as the time tends to $T_c$. This is problematic under quantization or measurement noise.
\end{itemize}

In the remainder of the chapter we develop a methodology to design controllers that converge to the origin with a predefined convergence time bound. Additionally, we provide sufficient conditions for our methodology to yield bounded time-varying gains. 

\section{Preliminaries and Problem Statement}
\label{sec:Prelim}
\subsection{Fixed-time stability and settling-time function}

Consider the system
\begin{equation}\label{eq:sys}
    \dot{\mathbf{x}}=\mathbf{f}(\mathbf{x},t)+\mathbf{b}_{n}d(t), \ \forall t\geq 0, 
\end{equation}
where $\mathbf{x}\in\mathbb{R}^n$ is the state of the system, $t\in[0,+\infty)$ is time, $\mathbf{b}_{n}=[0,\ldots,0,1]^T$, and $d$ is a disturbance satisfying $|d(t)|\leq \Delta$, for a constant $d<\infty$. \footnote{In the spirit of Filippov's interpretation of differential equations, solutions of~\eqref{eq:sys} are understood as any absolutely continuous function that satisfies the differential inclusion obtained by applying the Filippov regularization to $\mathbf{f}(\bullet, \bullet)$ (See \cite[Page 85]{Filippov1988DifferentialSides}), allowing us to consider $\mathbf{f}(\bullet, \bullet)$ discontinuous in the first argument. In the usual Filippov's interpretation, it is assumed that $\|\mathbf{f}(\mathbf{x}, t)\|$ has an integrable majorant function of time for any $\mathbf{x}$, ensuring existence and uniqueness of solutions in forward time. However, in this work we deal with $\mathbf{f}(\mathbf{x}, t)$ for which no majorant function exist, but existence and uniqueness of solutions is still guaranteed by argument similar to \cite{aldana2019design}. In particular, existence of solutions follows directly from the equivalence of solutions to a well-posed Filippov system via the time-scale transformation.}

The set of admissible disturbances is denoted by $\Pi$.
The solution of \eqref{eq:sys}, with disturbance $d$ and initial condition $\mathbf{x}_0$ is denoted by $\mathbf{x}(t;\mathbf{x}_0,d)$. If $d(t)\equiv0$ we simply write $\mathbf{x}(t;\mathbf{x}_0)$. Furthermore, consider the origin to be an equilibrium point of \eqref{eq:sys} for every admissible disturbance, meaning that $\mathbf{x}(t;\mathbf{0},d)=\mathbf{0}$ for all $t\ge 0$. 

\index{settling-time function}
\begin{definition}(Settling-time function)
\label{Def:Settling}
Then, the \textit{settling-time function} of system~\eqref{eq:sys} is defined as $T : \mathbb{R}^n \to \Bar{\mathbb{R}}_+$,
$$T(\mathbf{x}_0):=\inf\left\{\xi\geq 0:\forall d\in\Pi, \lim_{t\to\xi}\mathbf{x}(t;\mathbf{x}_0,d)=0\right\}.$$
\end{definition}
Notice that Definition~\ref{Def:Settling} admits $T(\mathbf{x}_0)=\infty$.

\index{Fixed-time stability}

\begin{definition} \label{def:finite}(Finite-time stability) 
The origin of system \eqref{eq:sys} is said to be \textit{finite-time stable} if it is asymptotically stable~\cite{Khalil2002NonlinearSystems} and its settling-time function is finite for every $\mathbf{x}_0$, i.e., $T(\mathbf{x}_0)<\infty$ for all $\mathbf{x}_0 \in \mathbb{R}^n$.
\end{definition}

\begin{definition} \label{def:fixed}(Fixed-time stability) 
The origin of system \eqref{eq:sys} is said to be \textit{fixed-time stable} if it is finite-time stable and its settling-time function $T(\mathbf{x}_0)$ is uniformly bounded on $\mathbb{R}^n$, i.e. there exists $T_{\text{max}}\in\mathbb{R}_+\setminus\{0\}$ such that $\sup_{\mathbf{x}_0\in\mathbb{R}^n}T(\mathbf{x}_0)\leq T_{\text{max}}$. Then, $T_{\text{max}}$ is said to be a \textit{UBST} of system \eqref{eq:sys}.
\end{definition}

\subsection{Problem Statement}
\label{sec:Problem}

Consider a chain of integrators 
\begin{align}
    \dot{x}_i&=x_{i+1}, \ \ \ \ i=1,\ldots,n-1 \\
    \dot{x}_n&=u(\mathbf{x},t;T_c)+d(t) \label{Eq:IntegratorChain}
\end{align}
where the disturbance $d(t)$ satisfies $|d(t)|\leq \Delta$ with a known constant $\Delta$,  $\mathbf{x}=[x_1,\ldots,x_n]^T$. We aim to design a controller $u(\mathbf{x},t;T_c)$ to steer the system to the origin before the desired time $T_c$ a priori set by the user, i.e., the controller $u(\mathbf{x},t;T_c)$ is such that the origin of the closed-loop system is fixed-time stable with a predefined \textit{UBST} given by $T_c$. 

\index{Predefined-time controller}
\index{Prescribed-time controller}
\begin{definition}
The controller $u(\mathbf{x},t;T_c)$ is called:
\begin{itemize}
    \item a \textit{predefined-time controller} if the settling-time function of the closed loop system satisfies $\sup_{\mathbf{x}_0\in\mathbb{R}^n}T(\mathbf{x}_0)\leq T_c<\infty$. 
    \item a \textit{prescribed-time controller} if for all $\mathbf{x}_0\neq0$ the settling-time function of the closed loop system satisfies $T(\mathbf{x}_0)=T_c<\infty$.
\end{itemize}
\end{definition}
Notice that prescribed-time controllers ensure convergence with an \textit{UBST} given by $T_c$. Thus, prescribed-time controllers are a subclass of predefined-time ones, with the remarkable property that the settling-time function is precisely $T_c$. 

\index{Hybrid controller for predefined-time convergence}
Our approach is a hybrid controller of the form:
\begin{equation}
\label{Eq:POC}
    u(\mathbf{x},t;T_c)=
    \left\lbrace 
    \begin{array}{lll}
     \phi(\mathbf{x},t;T_c)     &  \text{ for } & t\in[0,T_c)\\
      w(\mathbf{x};\Delta)  & & \text{otherwise}
    \end{array}
    \right.
\end{equation}
where the time-varying controller $\phi(\mathbf{x},t;T_c)$ should drive the state of the system to the origin with a convergence time bound specified a priori by the parameter $T_c$ and the robust controller $w(\mathbf{x};\Delta)$ should maintain the system at the origin in spite of the bounded disturbance $d(t)$. Since the design of robust sliding-mode controllers $w(\mathbf{x};\Delta)$ is well understood, see, e.g.,~\cite{Ding2016SimpleController,ShihongDing2015NewControllers,Utkin1992,Shtessel2014ObservationObservers}, in the rest of the chapter we focus on the design of the controller $\phi(\mathbf{x},t;T_c)$ and restrict the analysis to the interval $[0,T_c)$.

\index{Time-varying gain}
\section{First-order controllers}
\label{sec:first_order}
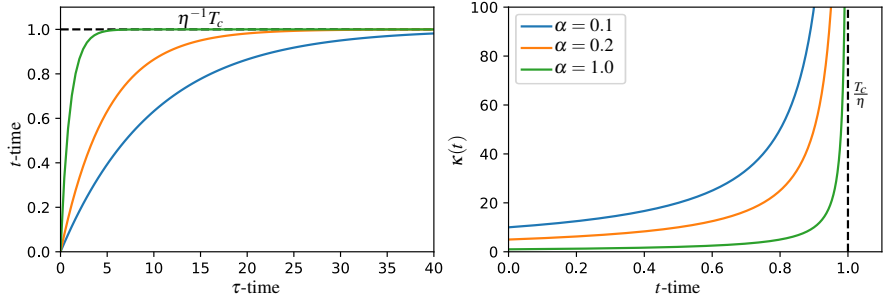
\begin{figure}
    \centering
\def\svgwidth{12cm}
\begingroup%
  \makeatletter%
  \providecommand\color[2][]{%
    \errmessage{(Inkscape) Color is used for the text in Inkscape, but the package 'color.sty' is not loaded}%
    \renewcommand\color[2][]{}%
  }%
  \providecommand\transparent[1]{%
    \errmessage{(Inkscape) Transparency is used (non-zero) for the text in Inkscape, but the package 'transparent.sty' is not loaded}%
    \renewcommand\transparent[1]{}%
  }%
  \providecommand\rotatebox[2]{#2}%
  \newcommand*\fsize{\dimexpr\f@size pt\relax}%
  \newcommand*\lineheight[1]{\fontsize{\fsize}{#1\fsize}\selectfont}%
  \ifx\svgwidth\undefined%
    \setlength{\unitlength}{613.203125bp}%
    \ifx\svgscale\undefined%
      \relax%
    \else%
      \setlength{\unitlength}{\unitlength * \real{\svgscale}}%
    \fi%
  \else%
    \setlength{\unitlength}{\svgwidth}%
  \fi%
  \global\let\svgwidth\undefined%
  \global\let\svgscale\undefined%
  \makeatother%
  \begin{picture}(1,0.35159206)%
    \lineheight{1}%
    \setlength\tabcolsep{0pt}%
    \put(0,0){\includegraphics[width=\unitlength]{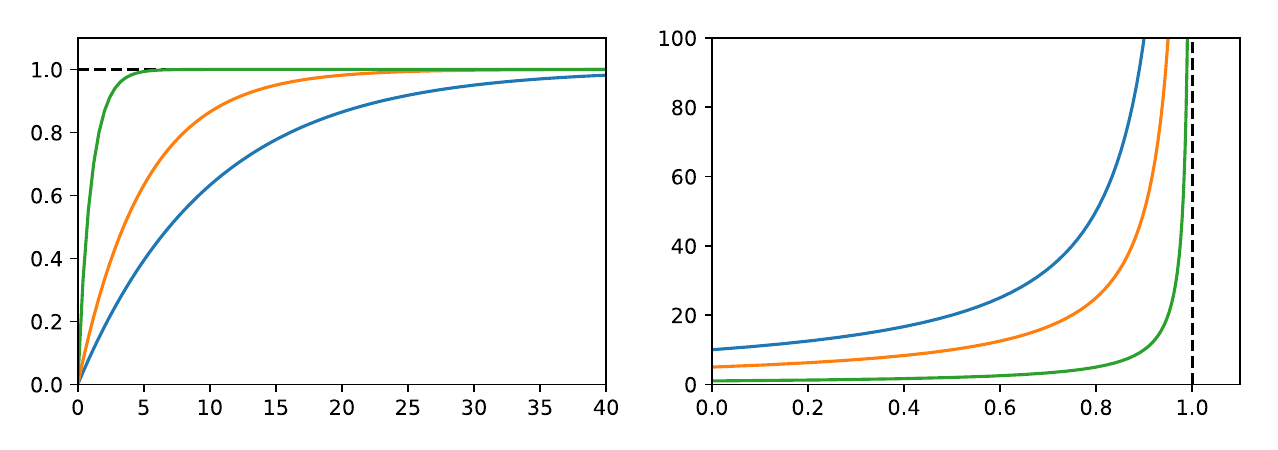}}%
    \put(0,0){\includegraphics[width=\unitlength]{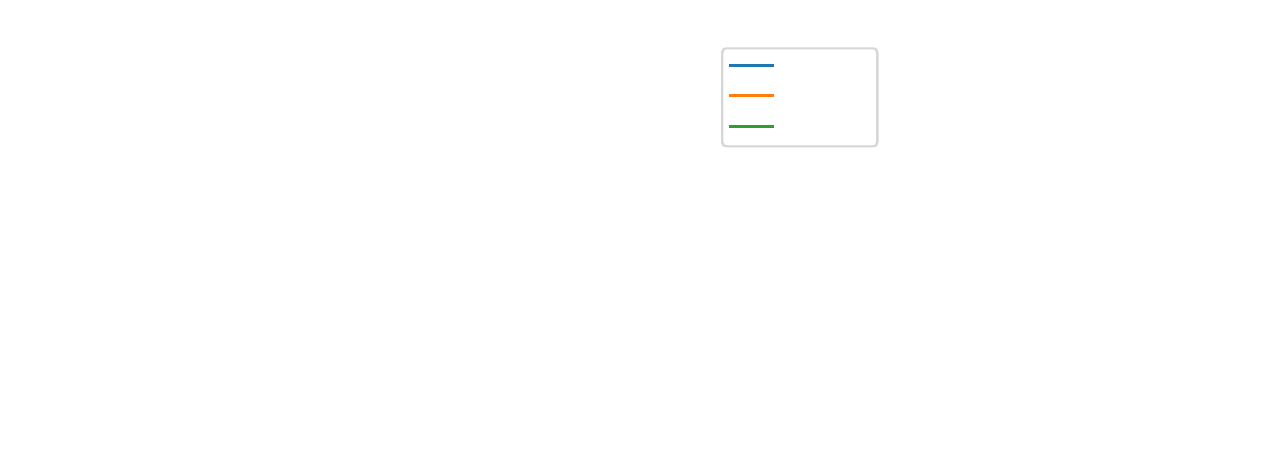}}%
    \scriptsize{
    \put(0.1902426,0.30322736){\makebox(0,0)[lt]{\lineheight{1.25}\smash{\begin{tabular}[t]{l}$\eta^{-1}T_c$\end{tabular}}}}%
    \put(0.25,0.00541022){\makebox(0,0)[lt]{\lineheight{1.25}\smash{\begin{tabular}[t]{l}$\tau$-time\end{tabular}}}}%
    \put(0.72,0.00541022){\makebox(0,0)[lt]{\lineheight{1.25}\smash{\begin{tabular}[t]{l}$t$-time\end{tabular}}}}%
    \put(0.0165864,0.14477868){\rotatebox{90}{\makebox(0,0)[lt]{\lineheight{1.25}\smash{\begin{tabular}[t]{l}$t$-time\end{tabular}}}}}%
    \put(0.50582071,0.14477868){\rotatebox{90}{\makebox(0,0)[lt]{\lineheight{1.25}\smash{\begin{tabular}[t]{l}$\kappa(t)$\end{tabular}}}}}%
    \put(0.60786844,0.29512914){\makebox(0,0)[lt]{\lineheight{1.25}\smash{\begin{tabular}[t]{l}$\alpha=0.1$\end{tabular}}}}%
    \put(0.60786844,0.2731136){\makebox(0,0)[lt]{\lineheight{1.25}\smash{\begin{tabular}[t]{l}$\alpha=0.2$\end{tabular}}}}%
    \put(0.60786844,0.24865189){\makebox(0,0)[lt]{\lineheight{1.25}\smash{\begin{tabular}[t]{l}$\alpha=1.0$\end{tabular}}}}%
    \put(0.93644364,0.22031615){\makebox(0,0)[lt]{\lineheight{1.25}\smash{\begin{tabular}[t]{l}$\frac{T_c}{\eta}$\end{tabular}}}}%
    }
  \end{picture}%
\endgroup%
    \caption{Example of a time-scale transformation (left) and its related time-varying gain (right) with $\eta=1$ and $T_c=1$.}
    \label{fig:TimeScale}
\end{figure}

\index{Time-scale transformation for predefined-time convergence}
Consider the time-scale transformation 
\begin{equation}
    \label{Eq:TimeScaleTransf}
    \tau=\varphi(t)=\ln((1-\eta T_c^{-1} t)^{-\frac{1}{\alpha}})
\end{equation}
with constant positive parameters $\alpha$, $\eta$, and $T_c$. Its inverse is given by 
$$t=\varphi^{-1}(\tau)=\eta^{-1}T_c\left(1-\exp(-\alpha \tau)\right),$$ together with the time-varying gain 
\begin{equation}
\label{Eq:kappa}
\kappa(t):=\frac{\mathrm{d}\tau}{\mathrm{d}t}=\frac{\eta}{\alpha(T_c-\eta t)}.
\end{equation}
Such a time-scale transformation is illustrated in Fig.~\ref{fig:TimeScale}. Notice that as $\tau$ tends to infinity $t$ approaches $\eta^{-1}T_c$, and $\lim_{t\to\eta^{-1}T_c^-}\kappa(t)=\infty$. This property will be exploited to design an asymptotically stable system in $\tau$-time and transform it into a predefined-time system in $t$-time, as explained next.

Consider a first-order system
\begin{equation}
    \label{Eq:FOS}
    \dot{x}=\phi(x,t;T_c)+d(t), \ \ \ \ \ t\in[0,T_c)
\end{equation}
with the controller
\begin{equation}
\label{Eq:FOPhi}
    \phi(x,t;T_c)=\kappa(t)v(x),
\end{equation}
where $x\in\mathbb{R}$, and $v(x)$ is a virtual control to be defined below. 
System~\eqref{Eq:FOS} in $\tau$-time is given by
\begin{equation}
    \frac{\mathrm{d}x}{\mathrm{d}\tau}=\left.\frac{\mathrm{d}x}{\mathrm{d}t}\right|_{t=\eta^{-1}T_c\left(1-\exp(-\alpha \tau)\right)}\cdot \frac{\mathrm{d}t}{\mathrm{d}\tau}
\end{equation}
where 
\begin{equation}
\label{Eq:dtdtau}
    \frac{\mathrm{d}t}{\mathrm{d}\tau}=\alpha\eta^{-1}T_c\exp(-\alpha\tau)=\kappa(\varphi^{-1}(\tau))^{-1}.
\end{equation}
Thus, 
\begin{equation}
\label{Eq:kappatau}
    \kappa(\varphi^{-1}(\tau))=\frac{\eta}{\alpha T_c}\exp(\alpha\tau)
\end{equation}
and
\begin{equation}
    \label{Eq:FOS-TauTime}
    \frac{\mathrm{d}x}{\mathrm{d}\tau}=v(x)+\alpha\eta^{-1}T_c\exp(-\alpha\tau)d(\varphi^{-1}(\tau)).
\end{equation}

Notice that, since $|d(\varphi^{-1}(\tau))|\leq \Delta$, then the disturbance term 
\begin{equation}
    \label{Eq:VanishDist}
\alpha\eta^{-1}T_c\exp(-\alpha\tau)d(\varphi^{-1}(\tau))
\end{equation}
 becomes vanishing in $\tau$-time. 
Thus, if $v(x)$ is chosen such that~\eqref{Eq:FOS-TauTime} is asymptotically stable with a settling-time function $\mathcal{T}(x_0)$, due to the time-scale transformation~\eqref{Eq:TimeScaleTransf}, the settling-time function of~\eqref{Eq:FOS} is given by
\begin{equation}
    \label{Eq:SttlingTimeMotiv}
    T(x_0)=\eta^{-1}T_c\left(1-\exp\left(-\alpha \mathcal{T}(x_0)\right)\right).
\end{equation}

Thus, by an appropriate selection of $v(x)$ and $\eta$, we can obtain a predefined-time controller~\eqref{Eq:POC}. 

One drawback of the controller~\eqref{Eq:FOPhi} is that, if $v(x)$ contains discontinuous terms, then $\phi(x,t;T_c)$  will have discontinuous terms that are increasing beyond what is necessary to cancel out the disturbance effect, possibly producing large chattering. For, instance, notice that with $v(x)=-\sign{x}$ and $\eta=1$,~\eqref{Eq:FOS-TauTime} is finite-time stable but $\phi(x,t;T_c)=-\frac{1}{\alpha(T_c-t)}\sign{x}$.

To address this important issue, consider the following generalization of the controller in~\eqref{Eq:FOPhi}, with an additional degree of freedom $\rho \in [0, 1]$:
\begin{equation}
\label{Eq:FOCrho}
    \phi(x,t;T_c)=
    \beta\kappa(t)^{1-\rho}\tilde{v}(\beta^{-1}\kappa(t)^\rho x)
\end{equation}
where $\beta\geq(\alpha\eta^{-1}T_c)^{1-\rho}$, $\kappa(t)$ is given in~\eqref{Eq:kappa}, and $\tilde{v}(\bullet)$ is an auxiliary controller to be specified below. 
To analyze the stability of the closed-loop system, consider the coordinate change:
$$
z=\beta^{-1}\kappa(t)^\rho x,
$$
together with the time-scale transformation in~\eqref{Eq:TimeScaleTransf}. Noticing that 
\begin{equation}
\label{Eq:k*kdot}
\dot{\kappa}(t)\kappa(t)^{-1}=\alpha\kappa(t),
\end{equation}
the dynamics in the $z$-coordinates is given by
\begin{align}
\dot{z}&=\rho \kappa(t)^{-1}\dot{\kappa}(t) z+ \kappa(t)\tilde{v}(z)+\beta^{-1}\kappa(t)^\rho d(t)\\
&=\kappa(t)\left(\rho \alpha z+ \tilde{v}(z)+\beta^{-1}\kappa(t)^{\rho-1} d(t)\right).
\end{align}

Thus, from~\eqref{Eq:dtdtau} and~\eqref{Eq:kappatau}, it follows that
the dynamics in $z$-coordinates and $\tau$-time is given by
\begin{equation}
    \label{Eq:FOS-TauTimerho}
    \frac{\mathrm{d}z}{\mathrm{d}\tau}=\rho\alpha z+\tilde{v}(z)+\beta^{-1}(\alpha\eta^{-1}T_c)^{1-\rho}\exp(-\alpha(1-\rho)\tau)d(\varphi^{-1}(\tau)).
\end{equation}
Notice that
$$
\pi(\tau)=\beta^{-1}(\alpha\eta^{-1}T_c)^{1-\rho}\exp(-\alpha(1-\rho)\tau)d(\varphi^{-1}(\tau))
$$
satisfies 
$$
|\pi(\tau)|\leq \Delta\exp(-\alpha(1-\rho)\tau)
$$
and, therefore, with the $\rho$ parameter, we can specify the rate at which $\pi(\tau)$ vanishes. Moreover, with $\rho=1$, $\pi(\tau)$ is no longer a vanishing disturbance.
Thus, choosing the auxiliary controller as: 
$$
\tilde{v}(z)=v(z)-\alpha\rho z
$$
yields
\begin{equation}
\label{Eq:FOSPi}
    \frac{\mathrm{d}z}{\mathrm{d}\tau}=v(z)+\pi(\tau).
\end{equation}
Thus, we can take advantage of existing robust controllers for~\eqref{Eq:FOSPi}, and the settling-time function will become
\begin{equation}
    \label{Eq:SttlingTimeMotivrho}
    T(x_0)=\eta^{-1}T_c\left(1-\exp\left(-\alpha \mathcal{T}(\beta^{-1}\kappa(0)^\rho x_0)\right)\right).
\end{equation}

Furthermore, with $\rho=1$, if $v(x)$ contains an additive discontinuous term (designed to cope with the disturbance $\pi(\tau)$), those terms will not be multiplied by $\kappa(t)$ in $\phi(x,t;T_c)$, and thus will not have its magnitude increased beyond what is necessary to reject the disturbance without increasing chattering. For instance, with $\rho=1$, $\eta=1$ and $v(x)=-\sign{x}$ we obtain:
$$
\phi(x,t;T_c)=
    -\beta\sign{x}-\frac{1}{
    (T_c-t)} x.
$$

\subsection{Prescribed-time controllers}
\label{Sec:Presc}
In this subsection, we focus on controllers $v(x)$, such that the settling-time function of the closed-loop system~\eqref{Eq:FOS-TauTime} satisfies $\mathcal{T}(x_0)=\infty, \forall x_0\neq 0$ and we choose $\eta=1$. Since in the $\tau$-time the disturbance becomes vanishing, then, any Input-to-State Stabilizing controller~\cite{Sontag2008InputResults} can be applied as $v(x)$ to stabilize system~\eqref{Eq:FOS-TauTime}, even without knowledge of the disturbance bound $\Delta$. This is because, for any bounded disturbance $d(\cdot)$, in $\tau$-time the disturbance term~\eqref{Eq:VanishDist} goes to zero as the $\tau$-time goes to infinity. However, knowledge on $\Delta$ is required to maintain the state at the origin after the time $T_c$.

Therefore, with the controller~\eqref{Eq:POC}, the settling-time of the closed-loop system~\eqref{Eq:FOS} is 
\begin{equation}
\label{Eq:STPresc}
    T(x_0)=T_c,
\end{equation}
i.e., the convergence occurs precisely at $T_c$ regardless of the initial condition $x_0$.

The following proposition provides a first-order prescribed-time controller with minimum-energy among all controllers driving the system state from $x(0)=x_0$ to $x(T_c)=0$.

\index{Minimum energy prescribed-time control}

\begin{proposition}
\label{Prop:MinimumEnergy}
Let $d(t)=0$. Then, the trajectory $x(t)$ resulting from controller~\eqref{Eq:POC} where  $v(x)=-x$, and $\kappa(t)$ is given in~\eqref{Eq:kappa}
 with $\alpha=1$ and $\eta=1$, under the constraints $x(0)=x_0$ and $x(T_c)=0$, minimizes the energy function 
\begin{equation}
\label{Eq:Energy}
     E_{T_c}=\int_0^{T_c}u(\xi)^2 \mathrm{d}\xi.
\end{equation} 
\end{proposition}
\begin{proof}
Using $\dot{x}(t)=u(t)$, one can build a Lagrangian for \eqref{Eq:Energy} as $L(t,x,\dot{x})=\dot{x}^2$. Hence, the well-known Euler-Lagrange equations~\cite[Page 38]{Liberzon2019CalculusTheory}: $$\frac{\mathrm{d}}{\mathrm{d}t}\left(\frac{\partial L}{\partial \dot{x}}\right)-\frac{\partial L}{\partial x} = 0$$
lead to $\ddot{x}=0$ or $\dot{x}=0, \forall t\in(0,T_c)$. Thus, the resulting trajectories which minimize \eqref{Eq:Energy} must be of the form $x(t)=c_1+c_2t, \forall t\in[0,T_c]$ for some constants $c_1,c_2$. This, along with boundary conditions $x(0)=x_0$, $x(T_c)=0$ leads to $x(t)=x_0(1-T_c^{-1}t)$, which satisfies $\dot{x}=\kappa(t)v(x)=-\frac{1}{T_c-t}x, \forall t\in[0,T_c)$, concluding the proof. 
\qed
\end{proof}

The main drawback of prescribed-time controllers is that the origin of~\eqref{Eq:FOS} is reached as the time-varying gain tends to infinity, which is problematic under noise or limited numerical precision. One may suggest, as a workaround to maintain the time-varying gain bounded, to consider, instead of controller~\eqref{Eq:POC}, the controller
\begin{equation}
\label{Eq:FOCWorkaround}
    u(x,t;T_c
    )=
    \left\lbrace 
    \begin{array}{lll}
    \phi(x,t;T_c)     &  \text{ for } & t\in[0,t_{\textrm{stop}})\\
      w(x;\Delta)  & & \text{otherwise}
    \end{array}
    \right.
\end{equation}
where $t_{\textrm{stop}}<T_c$. Unfortunately, with $\phi(x,t;T_c)=-\kappa(t)x$, the state at $t_{\textrm{stop}}$ grows linearly with $x_0$, as illustrated in the following example.

\begin{example}
Consider a prescribed-time controller with $v(x)=-x$, with $\eta=1$, $\alpha=1$ and $T_c=1$ and set $t_{\textrm{stop}}=0.9$. The trajectories for different initial conditions are shown in Fig.~\ref{Fig:LinearFOSKrstic}; notice that $x(t_{\textrm{stop}})=x_0(1-T_c^{-1}t_{\textrm{stop}})=0.1x_0$. A similar case occurs by taking \begin{equation}
\label{Eq:ControllerKamal}
    v(x)=c(1-\exp(-|x|))\sign{x},
\end{equation} 
where $c\geq 1$, with such controller, the origin of system~\eqref{Eq:FOS-TauTime} is asymptotically stable. Thus, we take $\eta=1$. A predefined-time controller~\eqref{Eq:POC} with $v(x)$ as in~\eqref{Eq:ControllerKamal} was proposed in~\cite{Pal2020DesignTime}.  The trajectories with $c=10$ for different initial conditions are shown in Fig.~\ref{Fig:KamalFOS}, in this case $x(t_{\textrm{stop}})$ is also an unbounded function of the initial condition $x_0$.

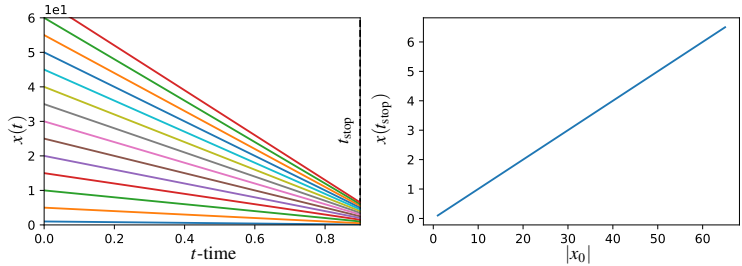
\begin{figure}
    \centering
\def\svgwidth{10cm}
\begingroup%
  \makeatletter%
  \providecommand\color[2][]{%
    \errmessage{(Inkscape) Color is used for the text in Inkscape, but the package 'color.sty' is not loaded}%
    \renewcommand\color[2][]{}%
  }%
  \providecommand\transparent[1]{%
    \errmessage{(Inkscape) Transparency is used (non-zero) for the text in Inkscape, but the package 'transparent.sty' is not loaded}%
    \renewcommand\transparent[1]{}%
  }%
  \providecommand\rotatebox[2]{#2}%
  \newcommand*\fsize{\dimexpr\f@size pt\relax}%
  \newcommand*\lineheight[1]{\fontsize{\fsize}{#1\fsize}\selectfont}%
  \ifx\svgwidth\undefined%
    \setlength{\unitlength}{603.66247559bp}%
    \ifx\svgscale\undefined%
      \relax%
    \else%
      \setlength{\unitlength}{\unitlength * \real{\svgscale}}%
    \fi%
  \else%
    \setlength{\unitlength}{\svgwidth}%
  \fi%
  \global\let\svgwidth\undefined%
  \global\let\svgscale\undefined%
  \makeatother%
  \begin{picture}(1,0.3684121)%
    \lineheight{1}%
    \setlength\tabcolsep{0pt}%
    \put(0,0){\includegraphics[width=\unitlength]{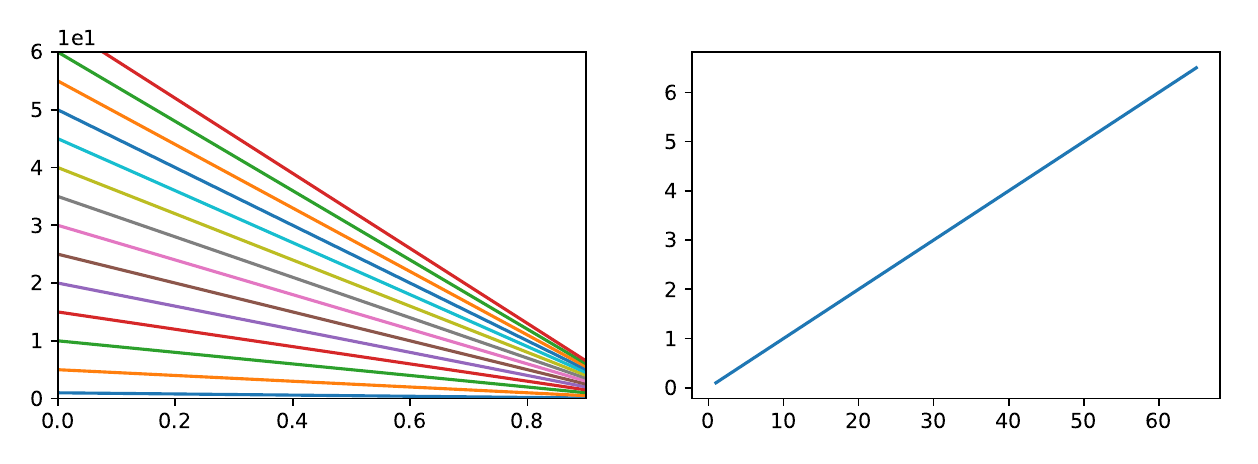}}%
    \scriptsize{
    \put(0.73856907,0.00623132){\color[rgb]{0,0,0}\makebox(0,0)[lt]{\lineheight{1.25}\smash{\begin{tabular}[t]{l}$|x_0|$\end{tabular}}}}%
    \put(0.24160262,0.00623132){\color[rgb]{0,0,0}\makebox(0,0)[lt]{\lineheight{1.25}\smash{\begin{tabular}[t]{l}$t$-time\end{tabular}}}}%
    \put(0.01693208,0.15959919){\color[rgb]{0,0,0}\rotatebox{90}{\makebox(0,0)[lt]{\lineheight{1.25}\smash{\begin{tabular}[t]{l}$x(t)$\end{tabular}}}}}%
    \put(0.45,0.15959919){\color[rgb]{0,0,0}\rotatebox{90}{\makebox(0,0)[lt]{\lineheight{1.25}\smash{\begin{tabular}[t]{l}$t_{\textrm{stop}}$\end{tabular}}}}}%
    \put(0.50,0.15959919){\color[rgb]{0,0,0}\rotatebox{90}{\makebox(0,0)[lt]{\lineheight{1.25}\smash{\begin{tabular}[t]{l}$x(t_{\textrm{stop}})$\end{tabular}}}}}%
    }
  \end{picture}%
\endgroup%
    \caption{Simulation of the first-order prescribed-time controller, for different initial conditions, with $\phi(x,t;T_c)=-\kappa(t)x$ and $T_c=1$. It can be observed that, the state $x(t_{\textrm{stop}})$ at a time $t_{\textrm{stop}}$ grows linearly with $|x_0|$. Here we choose $t_{\textrm{stop}}=0.9$.}
    \label{Fig:LinearFOSKrstic}
\end{figure}

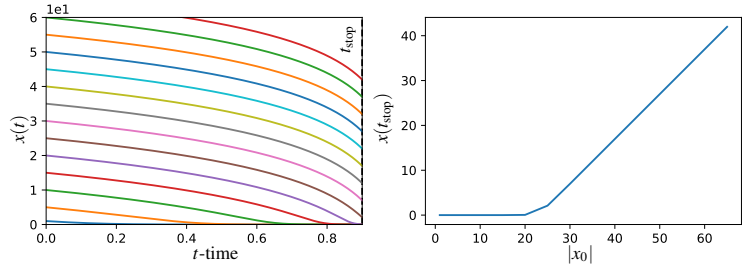
\begin{figure}
    \centering
\def\svgwidth{10cm}
\begingroup%
  \makeatletter%
  \providecommand\color[2][]{%
    \errmessage{(Inkscape) Color is used for the text in Inkscape, but the package 'color.sty' is not loaded}%
    \renewcommand\color[2][]{}%
  }%
  \providecommand\transparent[1]{%
    \errmessage{(Inkscape) Transparency is used (non-zero) for the text in Inkscape, but the package 'transparent.sty' is not loaded}%
    \renewcommand\transparent[1]{}%
  }%
  \providecommand\rotatebox[2]{#2}%
  \newcommand*\fsize{\dimexpr\f@size pt\relax}%
  \newcommand*\lineheight[1]{\fontsize{\fsize}{#1\fsize}\selectfont}%
  \ifx\svgwidth\undefined%
    \setlength{\unitlength}{603.66247559bp}%
    \ifx\svgscale\undefined%
      \relax%
    \else%
      \setlength{\unitlength}{\unitlength * \real{\svgscale}}%
    \fi%
  \else%
    \setlength{\unitlength}{\svgwidth}%
  \fi%
  \global\let\svgwidth\undefined%
  \global\let\svgscale\undefined%
  \makeatother%
  \begin{picture}(1,0.3684121)%
    \lineheight{1}%
    \setlength\tabcolsep{0pt}%
    \put(0,0){\includegraphics[width=\unitlength]{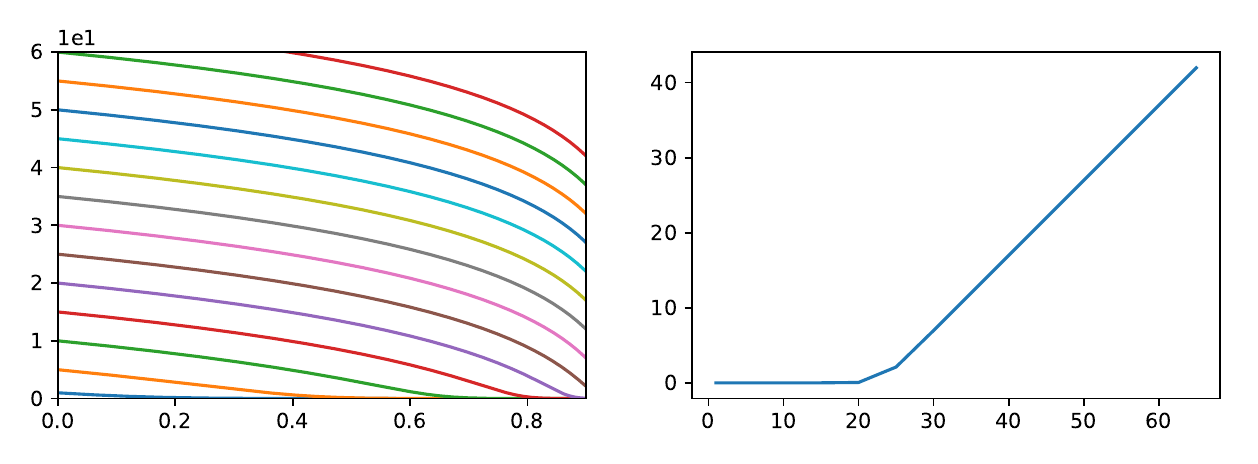}}%
    \scriptsize{
    \put(0.73856907,0.00623132){\color[rgb]{0,0,0}\makebox(0,0)[lt]{\lineheight{1.25}\smash{\begin{tabular}[t]{l}$|x_0|$\end{tabular}}}}%
    \put(0.24160262,0.00623132){\color[rgb]{0,0,0}\makebox(0,0)[lt]{\lineheight{1.25}\smash{\begin{tabular}[t]{l}$t$-time\end{tabular}}}}%
    \put(0.01693208,0.15959919){\color[rgb]{0,0,0}\rotatebox{90}{\makebox(0,0)[lt]{\lineheight{1.25}\smash{\begin{tabular}[t]{l}$x(t)$\end{tabular}}}}}%
    \put(0.45,0.28){\color[rgb]{0,0,0}\rotatebox{90}{\makebox(0,0)[lt]{\lineheight{1.25}\smash{\begin{tabular}[t]{l}$t_{\textrm{stop}}$\end{tabular}}}}}%
    \put(0.50,0.15959919){\color[rgb]{0,0,0}\rotatebox{90}{\makebox(0,0)[lt]{\lineheight{1.25}\smash{\begin{tabular}[t]{l}$x(t_{\textrm{stop}})$\end{tabular}}}}}%
    }
  \end{picture}%
\endgroup%
    \caption{Simulation of the first-order prescribed-time controller, for different initial conditions, with $\phi(x,t;T_c)=-\kappa(t)c(1-\exp(-|x|))\sign{x}$, with $c=10$ and $T_c=1$. It can be observed that, the state $x(t_{\textrm{stop}})$ at a time $t_{\textrm{stop}}$ grows with $|x_0|$. Here we choose $t_{\textrm{stop}}=0.9$.}
    \label{Fig:KamalFOS}
\end{figure}
\end{example}

\subsection{Predefined-time controllers with bounded time-varying gains}
\label{SubSec:PredefinedFOS}

As discussed above, prescribed-time controllers have the remarkable property that the settling-time function of the closed-loop system is precisely $T_c$. Still, they present a major drawback: the time-varying gain grows to infinity as the trajectory goes zero. Our approach to maintain the gain finite at the reaching time is to choose $v(x)$ such that $\mathcal{T}(x_0)<\infty$, i.e., such that the origin of~\eqref{Eq:FOS-TauTime} is finite-time stable. Then, the origin of~\eqref{Eq:FOS} is reached before the singularity in $\kappa(t)$ occurs. Moreover, a bounded time-varying gain can be obtained by choosing $v(x)$ such that $$\sup_{x_0\in\mathbb{R}}\mathcal{T}(x_0)\leq T_f<\infty$$
for a known $T_f$, i.e. such that the origin of~\eqref{Eq:FOS-TauTime} is fixed-time stable with a known \textit{UBST}.
Then, by choosing 
\begin{equation}
\label{Eq:Eta}
\eta=:\left(1-\exp\left(-\alpha T_f\right)\right)
\end{equation}
(notice that $\eta<1$) with the controller~\eqref{Eq:POC}, the origin is reached in a predefined-time $T_c$ and $\kappa(t)$ is bounded by 
\begin{equation}
    \label{Eq:kmax} 
    \kappa(t)\leq \kappa_{\textrm{max}}:=\frac{\exp(\alpha T_f)-1}{\alpha T_c} \ \ \text{ for } t\in[0,T_c),
\end{equation} 
with a settling-time function bounded by $T_c$. Moreover, if $\sup_{x_0\in\mathbb{R}}\mathcal{T}(x_0)= T_f$, then $\sup_{x_0\in\mathbb{R}}T(x_0)=T_c$.
Notice that the convergence is obtained before the desired time $T_c$, instead of precisely at time $T_c$, as with prescribed-time controllers.

\subsection{On reducing the energy with time-varying gains}

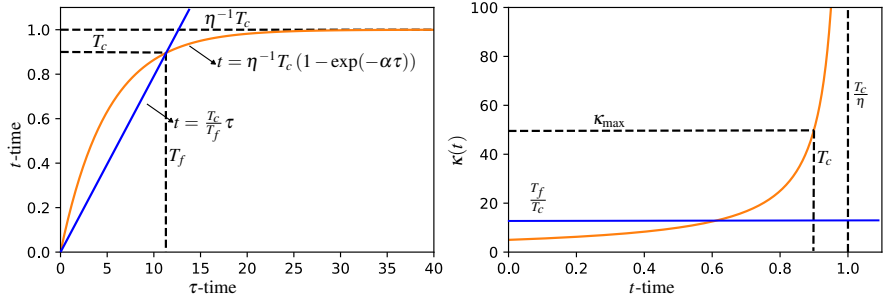
\begin{figure}
    \centering
\def\svgwidth{12cm}    
\begingroup%
  \makeatletter%
  \providecommand\color[2][]{%
    \errmessage{(Inkscape) Color is used for the text in Inkscape, but the package 'color.sty' is not loaded}%
    \renewcommand\color[2][]{}%
  }%
  \providecommand\transparent[1]{%
    \errmessage{(Inkscape) Transparency is used (non-zero) for the text in Inkscape, but the package 'transparent.sty' is not loaded}%
    \renewcommand\transparent[1]{}%
  }%
  \providecommand\rotatebox[2]{#2}%
  \newcommand*\fsize{\dimexpr\f@size pt\relax}%
  \newcommand*\lineheight[1]{\fontsize{\fsize}{#1\fsize}\selectfont}%
  \ifx\svgwidth\undefined%
    \setlength{\unitlength}{613.203125bp}%
    \ifx\svgscale\undefined%
      \relax%
    \else%
      \setlength{\unitlength}{\unitlength * \real{\svgscale}}%
    \fi%
  \else%
    \setlength{\unitlength}{\svgwidth}%
  \fi%
  \global\let\svgwidth\undefined%
  \global\let\svgscale\undefined%
  \makeatother%
  \begin{picture}(1,0.35159206)%
    \lineheight{1}%
    \setlength\tabcolsep{0pt}%
    \put(0,0){\includegraphics[width=\unitlength]{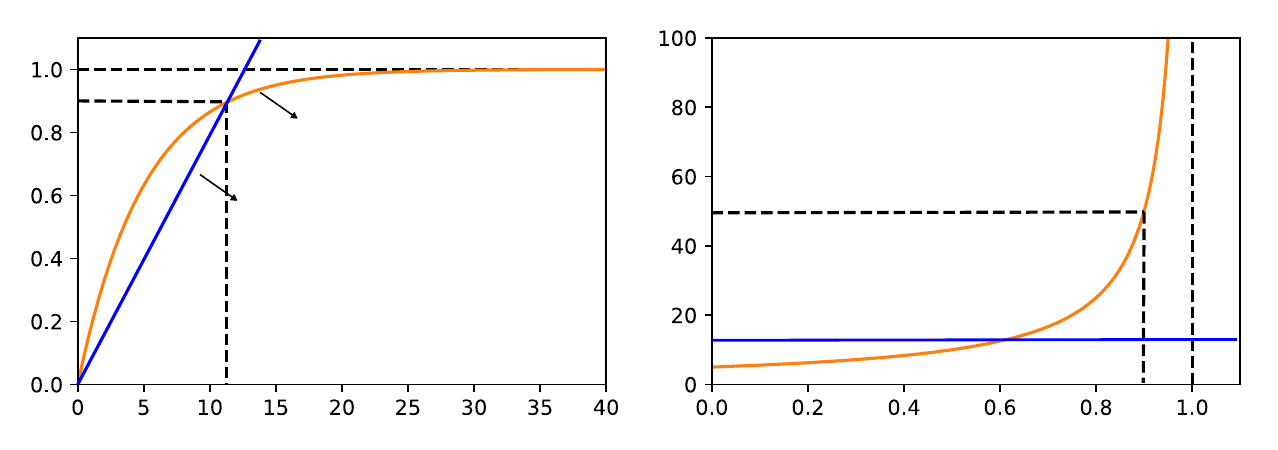}}%
    \scriptsize{
    \put(0.235,0.25582902){\makebox(0,0)[lt]{\lineheight{1.25}\smash{\begin{tabular}[t]{l}$t=\eta^{-1}T_c\left(1-\exp(-\alpha \tau)\right)$\end{tabular}}}}%
    \put(0.20218183,0.00541022){\makebox(0,0)[lt]{\lineheight{1.25}\smash{\begin{tabular}[t]{l}$\tau$-time\end{tabular}}}}%
    \put(0.6914161,0.00541022){\makebox(0,0)[lt]{\lineheight{1.25}\smash{\begin{tabular}[t]{l}$t$-time\end{tabular}}}}%
    \put(0.0165864,0.14477868){\rotatebox{90}{\makebox(0,0)[lt]{\lineheight{1.25}\smash{\begin{tabular}[t]{l}$t$-time\end{tabular}}}}}%
    \put(0.50582071,0.14477868){\rotatebox{90}{\makebox(0,0)[lt]{\lineheight{1.25}\smash{\begin{tabular}[t]{l}$\kappa(t)$\end{tabular}}}}}%
    \put(0.22,0.30322736){\makebox(0,0)[lt]{\lineheight{1.25}\smash{\begin{tabular}[t]{l}$\eta^{-1}T_c$\end{tabular}}}}%
    \put(0.09605545,0.27802234){\makebox(0,0)[lt]{\lineheight{1.25}\smash{\begin{tabular}[t]{l}$T_c$\end{tabular}}}}%
    \put(0.18095656,0.15199728){\makebox(0,0)[lt]{\lineheight{1.25}\smash{\begin{tabular}[t]{l}$T_f$\end{tabular}}}}%
    \put(0.89876308,0.15){\makebox(0,0)[lt]{\lineheight{1.25}\smash{\begin{tabular}[t]{l}$T_c$\end{tabular}}}}%
    \put(0.65002941,0.19250028){\makebox(0,0)[lt]{\lineheight{1.25}\smash{\begin{tabular}[t]{l}$\kappa_{\textrm{max}}$\end{tabular}}}}%
    \put(0.93644364,0.22031615){\makebox(0,0)[lt]{\lineheight{1.25}\smash{\begin{tabular}[t]{l}$\frac{T_c}{\eta}$\end{tabular}}}}%
    \put(0.18704769,0.18842276){\makebox(0,0)[lt]{\lineheight{1.25}\smash{\begin{tabular}[t]{l}$t=\frac{T_c}{T_f}\tau$\end{tabular}}}}%
    \put(0.57802566,0.10306473){\makebox(0,0)[lt]{\lineheight{1.25}\smash{\begin{tabular}[t]{l}$\frac{T_f}{T_c}$\end{tabular}}}}%
    }
  \end{picture}%
\endgroup%
    \caption{Comparison of the proposed time-scale transformation against the trivial scalar scaling. On the subplot on the left shows that if the system in $\tau$-time has an \textit{UBST} given by $T_f$, then the system in the $t$-time has a \textit{UBST} given by $T_c$. The subplot on the right illustrates how the time-varying gain is uniformly bounded.}
    \label{Fig:TVTimeScaleVsStaticTS}
\end{figure}

It is important to highlight that, in the first-order case, one can also obtain an autonomous predefined-time controller based on a fixed-time stable system with a known least \textit{UBST}, such as those proposed in~\cite{Sanchez-Torres2018,Aldana-Lopez2018}, by using the trivial time-scale transformation 
\begin{equation}
\label{Eq:TrivialTimeScaling}
    t=\frac{T_c}{T_f}\tau,
\end{equation} 
which result in the predefined-time controller $$u(t)=-\frac{T_f}{T_c}v(x)$$ 
where $T_c$ is the least \textit{UBST}. Fig.~\ref{Fig:TVTimeScaleVsStaticTS} illustrates how by using the time-scale transformation~\eqref{Eq:TimeScaleTransf}, the time-varying gain becomes bounded when a fixed-time controller $v(x)$ with \textit{UBST} given by $T_f$ is used, and it is contrasted with the static time-scale transformation~\eqref{Eq:TrivialTimeScaling} and its associated gain for predefined-time control.

Our approach yields this trivial time-scaling as a special case as $\alpha$ tends to zero, since $\lim_{\alpha\to 0}\varphi^{-1}(\tau)=\frac{T_c}{T_f}\tau$ and $\lim_{\alpha\to 0}\kappa(t)=\frac{T_f}{T_c}$.

As shown in the following proposition, even in the case where there already exists an autonomous fixed-time controller with least \textit{UBST}, our approach provides an extra degree of freedom to reduce the energy~\eqref{Eq:Energy},
used by the controller to drive system~\eqref{Eq:FOS} from $x(0)=x_0$ to $x(T_c)=0$, as well as to reduce the control magnitude $\sup_{t\in[0,T_c)}(|u(t)|)$. 

\begin{proposition}
Let the scalar system $\dot{x} = v(x)$ be such that its settling-time function satisfies $\sup_{x_0\in\mathbb{R}}\mathcal{T}(x_0)\leq T_f<\infty$
for a known $T_f$ and $v(\bullet)^2$ is non decreasing for non negative arguments. Using such $v(x), T_f$ and some $\alpha\geq 0$, construct a control $u(t)$ as in \eqref{Eq:POC} with $\eta$ defined in~\eqref{Eq:Eta} for the scalar system $\dot{x}=u$. Denote the energy $E(\alpha)=E_{T_c}$ as defined in \eqref{Eq:Energy} for such $\alpha\geq0$. Then, there always exist $\alpha^*_{x_0}\in\mathbb{R}$ which may depend on $x_0$ such that $E(\alpha^*_{x_0})\leq E(\alpha), \forall \alpha \geq 0$. In particular, $E(\alpha^*_{x_0})<E(0), \forall x_0\neq0$.
\end{proposition}
\begin{proof}
First, write the energy as $E(\alpha)=\int_0^{T_c}\kappa(\xi)^2v(x(\xi;x_0))^2d\xi$ using \eqref{Eq:Energy}, where $x(t;x_0)$ is the solution of $\dot{x} = u$ with $x(0;x_0)=x_0$. Now, make the change of variables $\tau = \varphi(\xi)$ from \eqref{Eq:TimeScaleTransf} which leads to 
$$
\begin{aligned}
E(\alpha) &= \int_0^{\varphi(T_c)} \kappa(\varphi^{-1}(\tau))^2v( x(\varphi^{-1}(\tau)) )^2 \left(\frac{1}{\kappa(\varphi^{-1}(\tau))}d\tau\right)\\
&=\int_0^{T_f}\left(\frac{\eta}{\alpha T_c}\exp(\alpha\tau)\right) v( x(\varphi^{-1}(\tau)) )^2 d\tau 
\end{aligned}
$$
where $d\tau = \kappa(\xi)d\xi$ was used from \eqref{Eq:kappa}. Now, note that $x(\varphi^{-1}(\tau))$ is the solution to $\frac{\mathrm{d}x}{\mathrm{d}\tau} = v(x)$ which follows from \eqref{Eq:FOS-TauTime} since there is no disturbance. Hence, $x(\varphi^{-1}(\tau))$ does not depend on $\alpha$. Now, it is straightforward to verify that:
$$
\lim_{\alpha\to0}\frac{\mathrm{d}}{\mathrm{d}\alpha} \left(\frac{\eta}{\alpha T_c}\exp(\alpha\tau)\right) = \frac{T_f}{T_c}(\tau - T_f/2)
$$
where we used $\eta=1-\exp\left(-\alpha T_f\right)$, from which it follows:
$$
\begin{aligned}
E'(0) &= \lim_{\alpha\to0}\frac{\mathrm{d}}{\mathrm{d}\alpha}\int_0^{T_f}\left(\frac{\eta}{\alpha T_c}\exp(\alpha\tau)\right)v(x(\tau))^2d\tau \\
&=\frac{T_f}{T_c}\int_0^{T_f}(\tau-T_f/2)v(x(\tau))^2d\tau\\
&=\frac{T_f}{T_c}\int_0^{T_f/2}(\tau-T_f/2)v(x(\tau))^2d\tau+\frac{T_f}{T_c}\int_{T_f/2}^{T_f}(\tau-T_f/2)v(x(\tau))^2d\tau \\
&=-\frac{T_f}{T_c}\int_0^{T_f/2}\tau v(x(T_f/2-\tau))^2d\tau+\frac{T_f}{T_c}\int_{0}^{T_f/2}\tau v(x(\tau+T_f/2))^2d\tau \\
&=\frac{T_f}{T_c}\int_0^{T_f/2}\tau \left(v(x(\tau+T_f/2))^2-v(x(T_f/2-\tau))^2\right)d\tau.
\end{aligned}
$$
Now, consider $x_0>0$. Hence, $x(\bullet)$ is a decreasing function and $x(\tau+T_f/2)\leq x(T_f/2-\tau)$ $\forall \tau\in[0,T_f/2]$. Therefore $v(x(\tau+T_f/2))^2\leq v(x(T_f/2-\tau))^2$ $\forall \tau\in[0,T_f/2]$. Hence, $E'(0)\leq 0$
with equality only if $x_0=0$ which is excluded in the proposition. Thus, one obtains the strict inequality $E(\alpha)<E(0)$ for $\alpha$ in some neighborhood of $0$. Now, note that $\lim_{\alpha\to\infty}\frac{\eta}{\alpha T_c}\exp(\alpha\tau)=+\infty$ so that $\lim_{\alpha\to\infty}E(\alpha)=+\infty$. Hence, there exists $\bar{\alpha}>0$ sufficiently big, such that $E(0)\leq E(\alpha), \forall \alpha\geq\bar{\alpha}$. Combining these facts, there must exist an optimal $\alpha^*>0$ such that $E(\alpha^*)\leq E(\alpha), \forall \alpha\geq 0$ due to continuity of $E(\alpha)$. The strict inequality $E(\alpha^*)<E(0)$ follows from $\alpha^*>0$ concluding the proof. 
\qed
\end{proof}

\begin{figure}
    \centering
\def\svgwidth{12cm}
\begingroup%
  \makeatletter%
  \providecommand\color[2][]{%
    \errmessage{(Inkscape) Color is used for the text in Inkscape, but the package 'color.sty' is not loaded}%
    \renewcommand\color[2][]{}%
  }%
  \providecommand\transparent[1]{%
    \errmessage{(Inkscape) Transparency is used (non-zero) for the text in Inkscape, but the package 'transparent.sty' is not loaded}%
    \renewcommand\transparent[1]{}%
  }%
  \providecommand\rotatebox[2]{#2}%
  \newcommand*\fsize{\dimexpr\f@size pt\relax}%
  \newcommand*\lineheight[1]{\fontsize{\fsize}{#1\fsize}\selectfont}%
  \ifx\svgwidth\undefined%
    \setlength{\unitlength}{613.3762207bp}%
    \ifx\svgscale\undefined%
      \relax%
    \else%
      \setlength{\unitlength}{\unitlength * \real{\svgscale}}%
    \fi%
  \else%
    \setlength{\unitlength}{\svgwidth}%
  \fi%
  \global\let\svgwidth\undefined%
  \global\let\svgscale\undefined%
  \makeatother%
  \begin{picture}(1,0.36257773)%
    \lineheight{1}%
    \setlength\tabcolsep{0pt}%
    \put(0,0){\includegraphics[width=\unitlength]{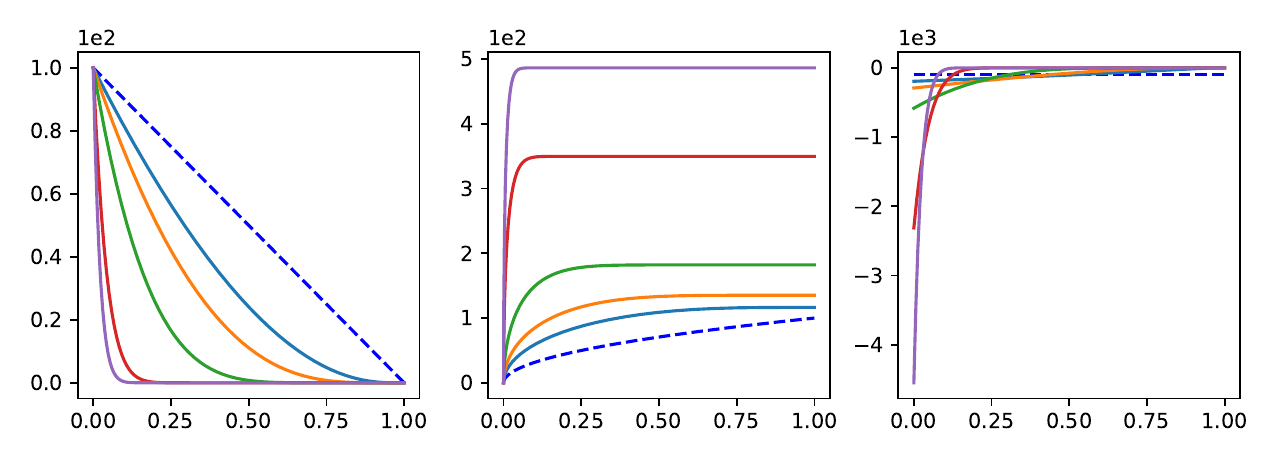}}%
    \put(0,0){\includegraphics[width=\unitlength]{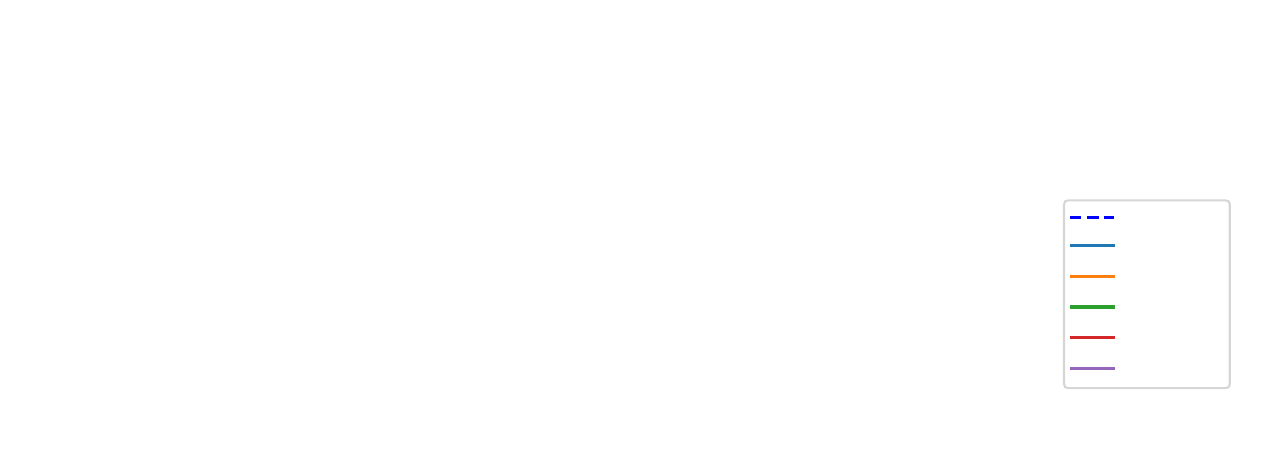}}%
    \scriptsize{
    \put(0.875,0.185){\makebox(0,0)[lt]{\lineheight{1.25}\smash{\begin{tabular}[t]{l}Optimal\end{tabular}}}}%
    \put(0.875,0.16602297){\makebox(0,0)[lt]{\lineheight{1.25}\smash{\begin{tabular}[t]{l}$\alpha=1.5$\end{tabular}}}}%
    \put(0.875,0.14156817){\makebox(0,0)[lt]{\lineheight{1.25}\smash{\begin{tabular}[t]{l}$\alpha=1.0$\end{tabular}}}}%
    \put(0.875,0.11711336){\makebox(0,0)[lt]{\lineheight{1.25}\smash{\begin{tabular}[t]{l}$\alpha=0.5$\end{tabular}}}}%
    \put(0.875,0.09265855){\makebox(0,0)[lt]{\lineheight{1.25}\smash{\begin{tabular}[t]{l}$\alpha=0.1$\end{tabular}}}}%
    \put(0.875,0.07064922){\makebox(0,0)[lt]{\lineheight{1.25}\smash{\begin{tabular}[t]{l}$\alpha=0.001$\end{tabular}}}}%
    \put(0.17090909,0.00475954){\makebox(0,0)[lt]{\lineheight{1.25}\smash{\begin{tabular}[t]{l}time\end{tabular}}}}%
    \put(0.48882165,0.00475954){\makebox(0,0)[lt]{\lineheight{1.25}\smash{\begin{tabular}[t]{l}time\end{tabular}}}}%
    \put(0.8067342,0.00475954){\makebox(0,0)[lt]{\lineheight{1.25}\smash{\begin{tabular}[t]{l}time\end{tabular}}}}%
    \put(0.0181179,0.15914732){\rotatebox{90}{\makebox(0,0)[lt]{\lineheight{1.25}\smash{\begin{tabular}[t]{l}$x(t)$\end{tabular}}}}}%
    \put(0.349,0.11){\rotatebox{90}{\makebox(0,0)[lt]{\lineheight{1.25}\smash{\begin{tabular}[t]{l}$E(t)=\left(\int_0^{t}|u(\xi)|^2 d\xi\right)^{\frac{1}{2}}$\end{tabular}}}}}%
    \put(0.668,0.15914732){\rotatebox{90}{\makebox(0,0)[lt]{\lineheight{1.25}\smash{\begin{tabular}[t]{l}$u(t)$\end{tabular}}}}}%
    }
  \end{picture}%
\endgroup%
    \caption{Simulation of a first-order predefined-time controller with bounded time-varying gains. Different values of the parameter alpha are shown, illustrating that a suitable selection of $\alpha$ allows to reduce the energy and the control magnitude.}
    \label{Fig:FOBilimit}
\end{figure}
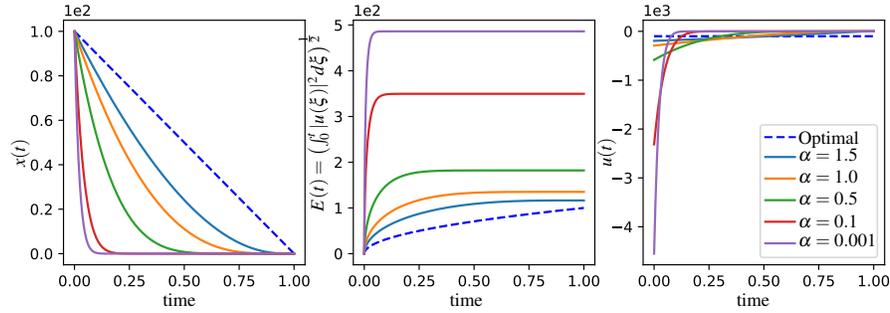

\begin{example}
\label{Ex:RedesignPolyakov}
Let $d(t)=0$ and the controller $v(x)=-\left((a|x|^{p} + b|x|^{q})^k\right) \sign{x}$, with $a$, $b$, $p$, $q$ and $k$ as in Theorem~\ref{th:tf_poly} in the appendix and $\zeta\geq\Delta$. Thus, $T_f=\gamma$ with $\gamma$ as in~\eqref{Eq:MinUpperEstimate} from that theorem. Then, with the controller~\eqref{Eq:POC} the origin is reached in a predefined-time $T_c$ and  $\kappa_{\textrm{max}}=\frac{\exp(\alpha \gamma)-1}{\alpha T_c}$. Fig.~\ref{Fig:FOBilimit} illustrates the trajectories for the case when $a=4$, $b=\frac{1}{4}$, $k=1$, $p=0.9$, $q=1.1$ (it follows from Theorem~\ref{th:tf_poly}, that $\gamma= 15.71$) obtained for different selections of $\alpha$, with $x_0=100$, as well as the Energy $E_{T_c}$ and the control signal $u(t)$ obtained in each case. Notice that, on the one hand, a minimum energy prescribed-time controller is obtained in Proposition~\ref{Prop:MinimumEnergy}, but requires time-varying gains that tend to infinity. On the other hand, a predefined-time autonomous controller can be obtained by taking $\lim_{\alpha\to 0}\kappa(t)$, but with an energy function significantly larger than with the prescribed-time controller. 
However, by tuning $\alpha$, the energy $ E_{T_c}$ and control magnitude $\sup_{t\in[0,T_c)}(|u(t)|)$ can be significantly reduced while maintaining the time-varying gain bounded. 
\end{example}

\subsection{Redesigning fixed-time stabilizing controllers using bounded time-varying gains: The second-order case}
\label{Sec:SOC}

In subsection~\ref{SubSec:PredefinedFOS} we argued that, even in the case where a fixed-time controller with a known least \textit{UBST} already exists, our approach allows to reduce the total energy required by the controller to drive the state of the system to the origin.

For higher-order systems, it is well known that the \textit{UBST} of fixed-time autonomous controllers, which are commonly based on Lyapunov analysis~\cite{Polyakov2012a} or homogeneity theory~\cite{Andrieu2008}, is very conservative, resulting in over-engineered controllers with large energy requirements
~\eqref{Eq:Energy} and large control signals. Thus, having predefined-time controllers that allow reducing such over-engineering is even more relevant in the high-order case.

Frequently, robust fixed-time controllers $v(\mathbf{z})$ are presented in the literature for system:
\begin{align}
\label{Eq:AuxSystem}
    \frac{\mathrm{d}z_1}{\mathrm{d}\tau}&=z_2\\
    \frac{\mathrm{d}z_2}{\mathrm{d}\tau}&= v(\mathbf{z})+d(t).
\end{align}
where $d(t)$ satisfy $|d|\leq\Delta$ for a known positive constant $\Delta$.
Assume that $v(\mathbf{z})$ is such that the origin of~\eqref{Eq:AuxSystem} is asymptotically stable and its settling-time function $\mathcal{T}(\mathbf{z}_0)$ satisfies  $$\sup_{\mathbf{z}_0\in\mathbb{R}^2}\mathcal{T}(\mathbf{z}_0)\leq T_f$$
for a known $T_f<\infty$. The approach described below can also be used to redesign finite-time controllers whose initial condition is bounded with a known settling-time function, see, e.g.,~\cite{Seeber2020ConvergenceControllers}.

To illustrate how to take advantage of such a controller and its estimation of the settling-time function to obtain predefined-time controllers based on time-varying gains, consider a more general predefined-time controller
\begin{equation}
\label{Eq:SOC}
    \phi(\mathbf{x},t;T_c)=\beta\kappa(t)^{2-\rho}\tilde{v}(\beta^{-1} \kappa(t)^{\rho}x_1,\beta^{-1} \kappa(t)^{\rho-1}x_2)
\end{equation}
where $\rho\in[0,2]$,  $\beta\geq (\alpha\eta^{-1}T_c)^{2-\rho}$, and $\tilde{v}(\bullet,\bullet)$ is an auxiliary control to be designed below.

To analyze the behavior of the closed loop system under controller~\eqref{Eq:SOC}, for $t\in[0,T_c)$, consider the coordinate change  
\begin{align}
\label{Eq:CoordXYSOS}
    y_i &= \beta^{-1} \kappa(t)^{\rho-i+1}x_i, \ \ \ i=1,2.
\end{align}

Recalling the equality~\eqref{Eq:k*kdot} then, the dynamics of the $\mathbf{y}$-coordinates are given by
\begin{align}
    \dot{y}_1
    &=\kappa(t)[\alpha\rho y_1+y_2]\\
    \dot{y}_2
    &=\kappa(t)[\alpha(\rho-1) y_2+\tilde{v}(y_1,y_2)+\beta^{-1} \kappa(t)^{\rho-2}d(t)].
\end{align}

Thus, considering the time-scale transformation~\eqref{Eq:TimeScaleTransf}, the dynamics in $\mathbf{y}$-coordinates and $\tau$-time are given by
\begin{align}
    \label{Eq:SOS-Tau}
    \frac{\mathrm{d}y_1}{\mathrm{d}\tau}&= \alpha\rho y_1+y_2\\
    \frac{\mathrm{d}y_2}{\mathrm{d}\tau}&= \alpha(\rho-1) y_2+\tilde{v}(y_1,y_2)+\pi(\tau)
\end{align}
where 
$$
\pi(\tau)=\beta^{-1} (\alpha\eta^{-1}T_c)^{2-\rho}\exp(-\alpha(2-\rho)\tau)d(\varphi^{-1}(\tau))
$$
satisfies $|\pi(\tau)|\leq\Delta$. Notice that, by choosing the $\rho$ parameter, the rate at which $\pi(\tau)$ vanishes can be varied; and with $\rho=2$, $\pi(\tau)$ becomes a non-vanishing disturbance.

Notice that, for the coordinate change~\eqref{Eq:CoordXYSOS} to be well defined, the auxiliary control $\tilde{v}(y_1,y_2)$, and the parameters $\rho$ and $\alpha$ need to be chosen such that 
\begin{equation}
    \lim_{\tau\to\infty} \exp(-\alpha(\rho+1-i)\tau) y_i(\tau;y_0,\pi_{[0,\infty]})=0 \ \ \ i=1,2,
\end{equation}
for every admissible disturbance $\pi_{[0,\infty]}$, which guarantees that the coordinate change maps the origin of the $\mathbf{y}$-coordinates to the origin of the $\mathbf{x}$-coordinates (and vice versa). Such condition is trivially satisfied since the origin of~\eqref{Eq:SOS-Tau} is finite-time stable.

Now, to design $\tilde{v}(\bullet,\bullet)$ based on the controller $v(\bullet)$ consider the coordinate change $z_1=y_1$ and $z_2=\alpha\rho y_1+y_2$, which will take the system into a controller canonical form:
\begin{align}
    \frac{\mathrm{d}z_1}{\mathrm{d}\tau}&=z_2\\
    \frac{\mathrm{d}z_2}{\mathrm{d}\tau}&=\alpha\rho z_2+\alpha(\rho-1) (z_2-\alpha\rho z_1)+\tilde{v}(z_1,z_2-\alpha\rho z_1)+\pi(\tau).
\end{align}
Thus, choosing $\tilde{v}(\bullet,\bullet)$ as 
$$
\tilde{v}(z_1,z_2-\alpha\rho z_1)=v(\mathbf{z})-c_1z_1-c_2z_2
$$
where
\begin{align}
\label{Eq:as}
    c_1&=-\alpha^2(\rho^2-\rho)\\
    c_2&=\alpha(2 \rho-1)
\end{align}
yields system~\eqref{Eq:AuxSystem}. Thus, taking $\eta$ as in~\eqref{Eq:Eta}, with the controller~\eqref{Eq:POC}
\begin{equation}
\label{Eq:SOC2}
\phi(\mathbf{x},t;T_c)=\beta\kappa(t)^{2-\rho}[v(\beta^{-1}\mathbf{Q}_{\rho} \mathbf{K}_{\rho}^{-1}(t) \mathbf{x})-[c_1,c_2]\beta^{-1}\mathbf{Q}_{\rho} \mathbf{K}_{\rho}^{-1}(t) \mathbf{x}]
\end{equation}
with $c_i$, $i=1,2$ is given by~\eqref{Eq:as}, $\mathbf{K}_{\rho}(t)=\mathrm{diag}( \kappa(t)^{-\rho}, \kappa(t)^{1-\rho})$
and
\begin{equation}
    \mathbf{Q}_{\rho}=
    \begin{bmatrix}
    1 & 0\\
    \alpha \rho & 1
    \end{bmatrix}
\end{equation}    
the origin of the closed-loop system is fixed-time stable and the settling-time function satisfies
\begin{equation}
    \label{Eq:Settling-timeSOS}
    T(\mathbf{x}_0)= \eta^{-1}T_c\left(1-\exp\left(-\alpha \mathcal{T}(\beta^{-1}\mathbf{Q}_{\rho} \mathbf{K}_{\rho}(0) \mathbf{x}_0)\right)\right)<T_c.
\end{equation}
Moreover, the time-varying gain is bounded by~\eqref{Eq:kmax}.

\begin{example}[Second-order system]
\label{Ex:SOSAut}

Consider the autonomous predefined-time controller given in Theorem~\ref{thm:socont}, where:
\begin{equation}
    \label{Eq:AutSOC}
u(t)=\omega(x_1,x_2)=-\left[\frac{\gamma_2}{T_{f_2}}\left(a_2\abs{\sigma}^{p}+b_2\abs{\sigma}^{q}\right)^{k}+\frac{\gamma_1^2}{2T_{f_1}^2}\left(a_1+3b_1x_1^2\right)+\zeta\right]\sign{\sigma}
\end{equation}
with
$$
\sigma=x_2+\barpow{\barpow{x_2}^2+\frac{2\gamma_1^2}{T_{f_1}^2}\left(a_1\barpow{x_1}^1+b_1\barpow{x_1}^3\right)}^{1/2},
$$
which was introduced in~\cite{Aldana-Lopez2018} for a second-order perturbed system. Consider a disturbance $d(t)=\sin(t)$, and let $a_1 = a_2 = 4$, $b_1 =b_2= \frac{1}{4}$, $p = 0.5$, $q = 1$, $k = 1.5$, $T_{f_1} =T_{f_2} = 5$, and $\zeta = 1$. According to Theorem~\ref{thm:socont} in the appendix, $\gamma_1$ and $\gamma_1$ are obtained as $\gamma_1=3.7081$ and $\gamma_2=2$, respectively, to obtain a predefined time controller with \textit{UBST} given by $T_c=T_{f_1}+T_{f_2}=10$. A simulation for initial conditions $x_1(0)=x_2(0)=50$ is shown in the first row of Fig.~\ref{Fig:ExSOS}.

Now, consider the predefined-time controller based on time-varying gains, given in~\eqref{Eq:SOC2}, using as a base controller $v(\mathbf{z})=\omega(z_1,z_2)$ the autonomous controller given in~\eqref{Eq:AutSOC}. Notice that $T_f=10$. 

Therefore,
$$
\mathbf{Q}_{\rho} \mathbf{K}_{\rho}^{-1}(t) \mathbf{x}=
\begin{bmatrix}
\kappa(t)^{\rho}x_1 \\
\alpha \rho \kappa(t)^{\rho}x_1+\kappa(t)^{\rho-1}x_2
\end{bmatrix}
$$
and

\begin{multline}
\phi(\mathbf{x},t;T_c)=\beta\kappa(t)^{2-\rho}\omega\left(\beta^{-1}\kappa(t)^{\rho}x_1,\beta^{-1}
\alpha \rho \kappa(t)^{\rho}x_1+\beta^{-1}\kappa(t)^{\rho-1}x_2\right)\\-(c 
_1+c_2
\alpha \rho )\kappa(t)^{2}x_1-c_2\kappa(t)x_2.
\end{multline}

We choose $\alpha=0.5$, $\rho=2$, $\beta=1$ and $T_c=10$. Thus, $\eta=0.9933$ and $\kappa_{\textrm{max}}=29.483$. 
Therefore,  the predefined-time controller~\eqref{Eq:POC} with $\phi(\mathbf{x},t;T_c)$ given by:

\begin{equation}
\phi(\mathbf{x},t;T_c)=\omega\left(\kappa(t)^{2}x_1,
 \kappa(t)^{2}x_1+\kappa(t)x_2\right)- \kappa(t)^{2}x_1-\frac{3}{2}\kappa(t)x_2
\end{equation}

A simulation for initial conditions $x_1(0)=x_2(0)=50$ is shown in the second row of Fig.~\ref{Fig:ExSOS}. Notice that, an improved transient behavior is obtained with the non-autonomous controller when compared with the behavior of the autonomous controller. Also notice that the control magnitude (second column) and the energy function (third column) are significantly reduced in the time-varying case.

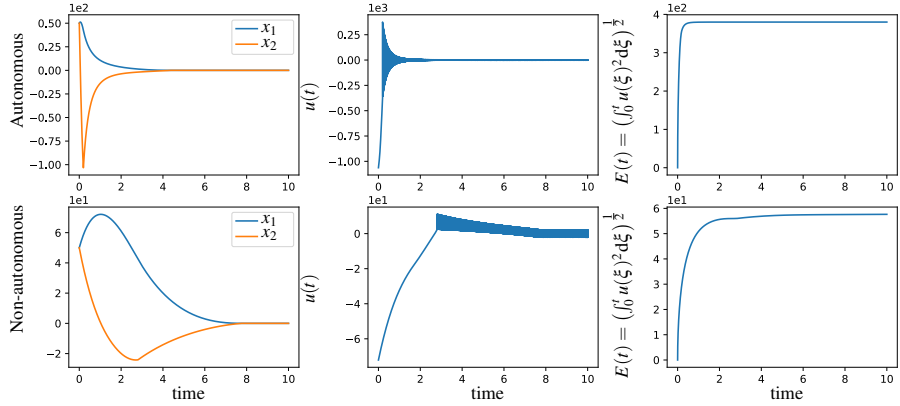
\begin{figure}
    \centering
\def\svgwidth{12cm}
\begingroup%
  \makeatletter%
  \providecommand\color[2][]{%
    \errmessage{(Inkscape) Color is used for the text in Inkscape, but the package 'color.sty' is not loaded}%
    \renewcommand\color[2][]{}%
  }%
  \providecommand\transparent[1]{%
    \errmessage{(Inkscape) Transparency is used (non-zero) for the text in Inkscape, but the package 'transparent.sty' is not loaded}%
    \renewcommand\transparent[1]{}%
  }%
  \providecommand\rotatebox[2]{#2}%
  \newcommand*\fsize{\dimexpr\f@size pt\relax}%
  \newcommand*\lineheight[1]{\fontsize{\fsize}{#1\fsize}\selectfont}%
  \ifx\svgwidth\undefined%
    \setlength{\unitlength}{851.14532471bp}%
    \ifx\svgscale\undefined%
      \relax%
    \else%
      \setlength{\unitlength}{\unitlength * \real{\svgscale}}%
    \fi%
  \else%
    \setlength{\unitlength}{\svgwidth}%
  \fi%
  \global\let\svgwidth\undefined%
  \global\let\svgscale\undefined%
  \makeatother%
  \begin{picture}(1,0.45669822)%
    \lineheight{1}%
    \setlength\tabcolsep{0pt}%
    \put(0,0){\includegraphics[width=\unitlength]{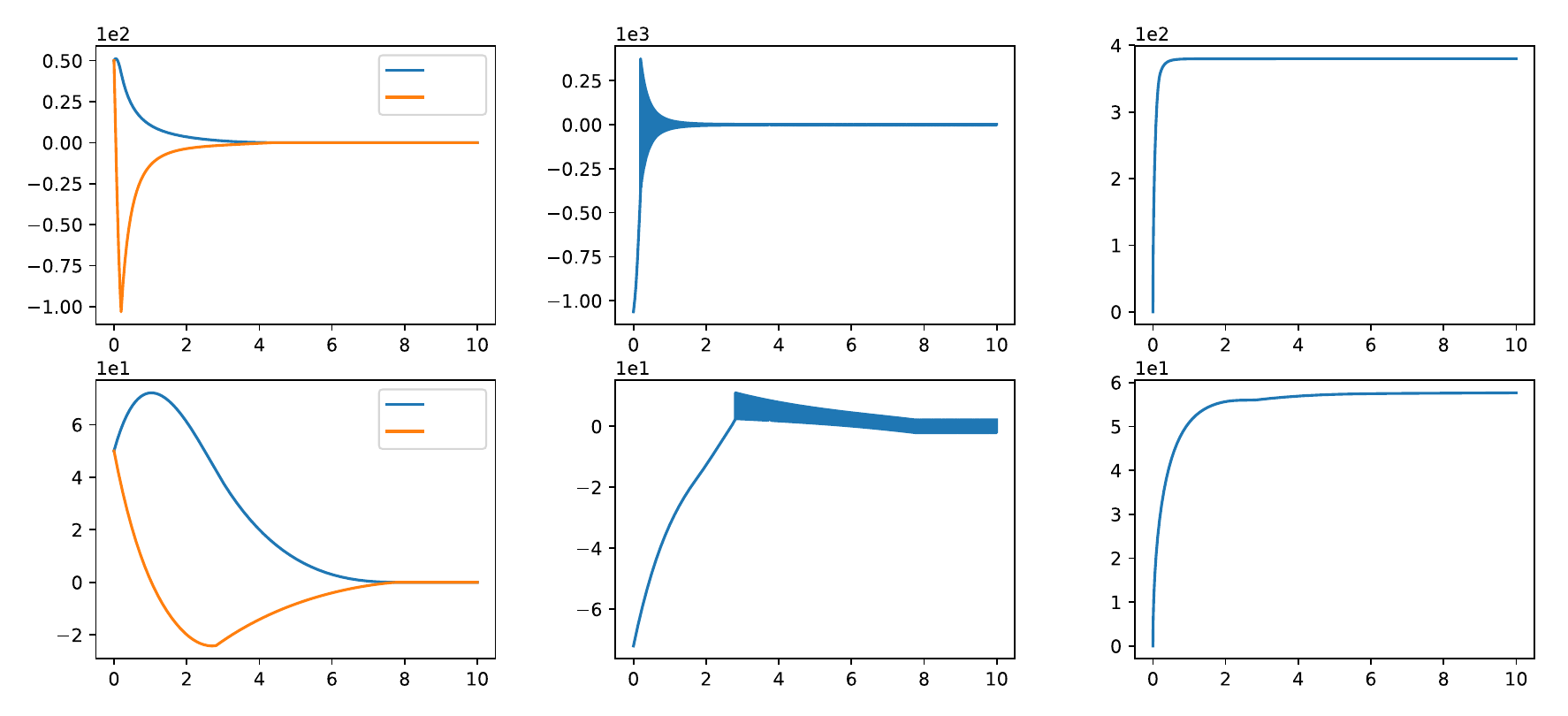}}%
    \scriptsize{
    \put(0.17550561,0.00106431){\makebox(0,0)[lt]{\lineheight{1.25}\smash{\begin{tabular}[t]{l}time\end{tabular}}}}%
    \put(0.50694137,0.00106431){\makebox(0,0)[lt]{\lineheight{1.25}\smash{\begin{tabular}[t]{l}time\end{tabular}}}}%
    \put(0.83837717,0.00106431){\makebox(0,0)[lt]{\lineheight{1.25}\smash{\begin{tabular}[t]{l}time\end{tabular}}}}%
    \put(0.01015928,0.30169082){\rotatebox{90}{\makebox(0,0)[lt]{\lineheight{1.25}\smash{\begin{tabular}[t]{l}Autonomous\end{tabular}}}}}%
    \put(0.01015928,0.07457606){\rotatebox{90}{\makebox(0,0)[lt]{\lineheight{1.25}\smash{\begin{tabular}[t]{l}Non-autonomous\end{tabular}}}}}%
    \put(0.33332614,0.10724968){\rotatebox{90}{\makebox(0,0)[lt]{\lineheight{1.25}\smash{\begin{tabular}[t]{l}$u(t)$\end{tabular}}}}}%
    \put(0.33332614,0.32042129){\rotatebox{90}{\makebox(0,0)[lt]{\lineheight{1.25}\smash{\begin{tabular}[t]{l}$u(t)$\end{tabular}}}}}%
    \put(0.68050539,0.025){\rotatebox{90}{\makebox(0,0)[lt]{\lineheight{1.25}\smash{\begin{tabular}[t]{l}$E(t)=\left(\int_0^t u(\xi)^2 \mathrm{d}\xi\right)^{\frac{1}{2}}$\end{tabular}}}}}%
    \put(0.68050539,0.24){\rotatebox{90}{\makebox(0,0)[lt]{\lineheight{1.25}\smash{\begin{tabular}[t]{l}$E(t)=\left(\int_0^t u(\xi)^2 \mathrm{d}\xi\right)^{\frac{1}{2}}$\end{tabular}}}}}%
    \put(0.275,0.40778614){\makebox(0,0)[lt]{\lineheight{1.25}\smash{\begin{tabular}[t]{l}$x_1$\end{tabular}}}}%
    \put(0.275,0.39192516){\makebox(0,0)[lt]{\lineheight{1.25}\smash{\begin{tabular}[t]{l}$x_2$\end{tabular}}}}%
    \put(0.275,0.19630639){\makebox(0,0)[lt]{\lineheight{1.25}\smash{\begin{tabular}[t]{l}$x_1$\end{tabular}}}}%
    \put(0.275,0.18044541){\makebox(0,0)[lt]{\lineheight{1.25}\smash{\begin{tabular}[t]{l}$x_2$\end{tabular}}}}%
    }
  \end{picture}%
\endgroup%
    \caption{Comparison between the autonomous predefined-time controller proposed in~\cite{Aldana-Lopez2018} and the proposed non-autonomous predefined-time controller as discussed in Example~\ref{Ex:SOSAut}. In both cases the \textit{UBST} is selected as $T_c=10$.}
    \label{Fig:ExSOS}
\end{figure}
\end{example}

The tuning parameters of our redesign methodology are $\alpha$, $\beta$, and $\rho$. The insight on the function of such parameters on the redesigned controller is as follows: The $\alpha$ parameter is associated with the time-scale transformation as illustrated in Fig.~\ref{fig:TimeScale}, increasing its value helps to reduce the slack between the \textit{UBST} and the true settling time~\cite{Gomez2020RNC}, however, it increases $\kappa_{\max}$, which impacts on chattering; Increasing the $\beta$ parameter allows to cope with disturbances $d(t)$ of greater magnitude, as in the $\tau$-time the magnitude of the disturbance $\pi(\tau)$ is inversely proportional to the magnitude of $\beta$.
However, large values of $\beta$ increase chattering; Finally, as mentioned above, the $\rho$ parameter allows to reduce chattering of the resulting predefined-time controller, in particular, when $\rho=n$ discontinuous terms in the admissible auxiliary controller, does not appear multiplied by the increasing time-varying gain in the redesigned controller.

\section{Main Result: Arbitrary-order predefined-time controller}
\label{sec:Main}

In this section, we present the extension to design arbitrary order predefined-time controllers. Our approach, can be seen as a ``redesign" methodology that starts from an admissible auxiliary controller, and uses time-varying gains to achieve predefined-time convergence. Let us introduce the notion of an admissible auxiliary controller.

\begin{definition}[admissible auxiliary controller]
Given parameters $\alpha\geq 0$, $\rho\in[0,n]$ and $T_f\in\bar{\mathbb{R}}_+$, we say that $v(\mathbf{z})$ is an \textit{admissible auxiliary controller} if:
\begin{enumerate}
    \item[(i)] for every disturbance $\pi(\tau)$ such that 
\begin{equation}
    \label{Eq:PiAux}
|\pi(\tau)|\leq\Delta\exp(-\alpha(n-\rho)\tau),
\end{equation}
it happens that the system
\begin{align}
\label{Eq:AuxSystemAOS}
\frac{\mathrm{d}\mathbf{z}}{\mathrm{d}\tau}
&=\mathbf{J}\mathbf{z}+\mathbf{b}_{n}(v(\mathbf{z})+\pi(\tau))
\end{align}
is asymptotically stable, where $\mathbf{J}:=[a_{ij}]\in\mathbb{R}^{n \times n}$ denotes a single Jordan block with zero eigenvalue, i.e., a square matrix with $a_{ij}=1$ if $j=i+1$ and $a_{ij}=0$ otherwise, and the vector $\mathbf{b}_i\in\mathbb{R}^n$ denotes a vector with one in the $i$-th entry and zeros otherwise. Moreover, the settling-time function satisfies
\begin{equation}
\label{Eq:STFAOS}
    \sup_{\mathbf{z}_0\in\mathbb{R}^n}\mathcal{T}(\mathbf{z}_0)\leq T_f.
\end{equation} 
\item[(ii)] for every admissible disturbance $\pi_{[0,\infty]}\in\Pi_{[0,\infty]}$, it happens that
\begin{equation}
\label{Eq:Limit}
    \lim_{\tau\to\infty} \exp(-\alpha(\rho+1-i)\tau) z_i(\tau;\mathbf{z}_0,\pi_{[0,\infty]})=0 \ \ \ i=1,\ldots,n.
\end{equation}
\end{enumerate}
\end{definition}

Based on an admissible auxiliary controller, we next present the methodology to design predefined-time controllers. Our main result is as follows.

\begin{theorem}
\label{Th:Main}
Given parameters $\alpha\geq 0$, $\rho\in[0,n]$ and $T_f\in\bar{\mathbb{R}}_+$, an admissible auxiliary controller $v(\mathbf{z})$, and a desired convergence time $T_c>0$, define the matrices $\mathbf{D}_{\rho},\mathbf{Q}_{\rho}\in\mathbb{R}^{n\times n}$ as: 
\begin{equation}
\label{Eq:MatrixDp}
\mathbf{D}_{\rho}:=\text{diag}\{-\rho,1-\rho,\ldots,n-1-\rho\}
\end{equation} 
and
\begin{equation}
\label{Eq:MatrixQp}
\mathbf{Q}_{\rho}: = \begin{bmatrix}
\mathbf{b}_{1}^{T}\\
\mathbf{b}_{1}^{T}(\mathbf{J}-\alpha \mathbf{D}_{\rho})\\
\vdots \\ 
\mathbf{b}_{1}^{T}(\mathbf{J}-\alpha \mathbf{D}_{\rho})^{n-1}
\end{bmatrix},
\end{equation}
and the time-varying matrix $\mathbf{K}_{\rho}(t)$ as
\begin{equation}
\label{Eq:Lambda}
\mathbf{K}_{\rho}(t):=\textrm{diag}(\kappa(t)^{-\rho},\kappa(t)^{1-\rho},\ldots,\kappa(t)^{n-1-\rho}),
\end{equation}
where $\kappa(t)$ is given by~\eqref{Eq:kappa} with $\eta$ as defined in~\eqref{Eq:Eta}. Then, with $\phi(\mathbf{x},t;T_c)$ given by:
\begin{equation}
\label{Eq:AOCPhi}
    \phi(\mathbf{x},t;T_c)=
    \beta\kappa(t)^{n-\rho}[v(\beta^{-1}\mathbf{Q}_{\rho} \mathbf{K}_{\rho}^{-1}(t) \mathbf{x})-\beta^{-1}\mathbf{b}_{1}^{T}(\mathbf{J}-\alpha \mathbf{D}_{\rho})^{n} \mathbf{K}_{\rho}^{-1}(t) \mathbf{x}]
\end{equation}
where $\beta\geq (\alpha\eta^{-1}T_c)^{n-\rho}$, the hybrid controller in ~\eqref{Eq:POC} is fixed-time stable with a settling time function given by
\begin{equation}
    \label{Eq:Settling-timeAOS}
    T(\mathbf{x}_0)= \eta^{-1}T_c\left(1-\exp\left(-\alpha \mathcal{T}(\beta^{-1}\mathbf{Q}_{\rho} \mathbf{K}_{\rho}(0) \mathbf{x}_0)\right)\right).
\end{equation}
\end{theorem}

\begin{proof}
Our approach for the proof of Theorem~\ref{Th:Main} is to show that the auxiliary system~\eqref{Eq:AuxSystemAOS} and the closed-loop system~\eqref{Eq:IntegratorChain} under controller~\eqref{Eq:POC},
in the time interval $[0,T_c)$, are related by the coordinate change
\begin{equation}
\label{eq:transformation}
\mathbf{z}=\beta^{-1}\mathbf{Q}_{\rho} \mathbf{K}_{\rho}^{-1}(t) \mathbf{x}
\end{equation}
and the time-scale transformation~\eqref{Eq:TimeScaleTransf}. 

Consider the time interval $[0,T_c)$ and the time-varying coordinate change~\eqref{eq:transformation}, and notice that, since $v(\mathbf{z})$ is an admissible auxiliary controller, then the coordinate change is well defined.
Then, the dynamics in the $\mathbf{z}$-coordinates is given by
\begin{equation}
\dot{\mathbf{z}}=\beta^{-1}\mathbf{Q}_{\rho} \frac{\mathrm{d}\mathbf{K}_{\rho}^{-1}(t)}{\mathrm{d}t} \mathbf{x}+\beta^{-1}\mathbf{Q}_{\rho} \mathbf{K}_{\rho}(t) \dot{\mathbf{x}}.
\end{equation}

Thus, it follows from the identity~\eqref{Eq:kID1} and $\mathbf{x}=\beta\mathbf{K}_{\rho}(t)\mathbf{Q}_{\rho}^{-1} \mathbf{z}$, that
\begin{align}
\dot{\mathbf{z}}&=-\alpha \kappa(t)\mathbf{Q}_{\rho}  \mathbf{D}_{\rho}\mathbf{Q}_{\rho}^{-1} \mathbf{z}+\beta^{-1}\mathbf{Q}_{\rho} \mathbf{K}_{\rho}(t) [\mathbf{J}\mathbf{x}+\mathbf{b}_{n}(u+d(t))]\\
&=-\alpha \kappa(t)\mathbf{Q}_{\rho}  \mathbf{D}_{\rho}\mathbf{Q}_{\rho}^{-1} \mathbf{z}+ \mathbf{Q}_{\rho} \mathbf{K}_{\rho}^{-1}(t)\mathbf{J}\mathbf{K}_{\rho}(t)\mathbf{Q}_{\rho}^{-1}\mathbf{z}+\beta^{-1}\mathbf{Q}_{\rho} \mathbf{K}_{\rho}^{-1}(t)\mathbf{b}_{n}(u+d(t)).
\end{align}
Moreover, applying identities~\eqref{Eq:QpId1} and~\eqref{Eq:kID2} from Lemma~\ref{Lem:DQ:identities} and Lemma~\ref{Lem:kappa:identities}, yields
\begin{align}
\dot{\mathbf{z}}&=-\alpha \kappa(t)\mathbf{Q}_{\rho}  \mathbf{D}_{\rho}\mathbf{Q}_{\rho}^{-1} \mathbf{z}+ \kappa(t)\mathbf{Q}_{\rho} \mathbf{J}\mathbf{Q}_{\rho}^{-1}\mathbf{z}+\beta^{-1}\mathbf{Q}_{\rho} \mathbf{K}_{\rho}^{-1}(t)\mathbf{b}_{n}(u+d(t))\\
&=\kappa(t)[ \mathbf{Q}_{\rho}  (\mathbf{J}-\alpha\mathbf{D}_{\rho})\mathbf{Q}_{\rho}^{-1} \mathbf{z}]+\beta^{-1}\mathbf{Q}_{\rho} \mathbf{K}_{\rho}^{-1}(t)\mathbf{b}_{n}(u+d(t))\\
&=\kappa(t)(\mathbf{J}+\mathbf{A})\mathbf{z}+\beta^{-1}\mathbf{Q}_{\rho} \mathbf{K}_{\rho}^{-1}(t)\mathbf{b}_{n}(u+d(t)).
\end{align}
From the identity $\mathbf{K}_{\rho}^{-1}(t)\mathbf{b}_{n}=\kappa(t)^{\rho-n+1}\mathbf{b}_{n}$ and $\mathbf{Q}_{\rho} \mathbf{b}_{n}=\mathbf{b}_{n}$ it follows that
\begin{align}
\dot{\mathbf{z}}&=\kappa(t)(\mathbf{J}+\mathbf{A})\mathbf{z}+\kappa(t)^{\rho-n+1}\beta^{-1}\mathbf{Q}_{\rho} \mathbf{b}_{n}(u+d(t))\\
&=\kappa(t)(\mathbf{J}+\mathbf{A})\mathbf{z}+\kappa(t)^{\rho-n+1}\beta^{-1}\mathbf{b}_{n}(u+d(t))\\
&=\kappa(t)(\mathbf{J}+\mathbf{A})\mathbf{z}+\kappa(t)\mathbf{b}_{n}v(\mathbf{z})-\mathbf{b}_{n}\mathbf{b}_{1}^{T}(\mathbf{J}-\alpha \mathbf{D}_{\rho})^{n}\mathbf{Q}_{\rho}^{-1} \mathbf{z}+\kappa(t)^{\rho-n+1}\beta^{-1}\mathbf{b}_{n}d(t)).
\end{align}
Since, according to~\eqref{Eq:A} from Lemma~\ref{Lem:DQ:identities}, $\mathbf{A}=\mathbf{b}_{n}\mathbf{b}_{1}^{T}(\mathbf{J}-\alpha \mathbf{D}_{\rho})^{n}\mathbf{Q}_{\rho}^{-1}$, then 
\begin{align}
\dot{\mathbf{z}}&=\kappa(t)(\mathbf{J}+\mathbf{A})\mathbf{z}+\kappa(t)\mathbf{b}_{n}v(\mathbf{z})-\mathbf{A} \mathbf{z}+\kappa(t)^{\rho-n+1}\beta^{-1}\mathbf{b}_{n}d(t))\\
&=\kappa(t)[\mathbf{J}\mathbf{z}+\mathbf{b}_{n}v(\mathbf{z})+\kappa(t)^{\rho-n}\beta^{-1}\mathbf{b}_{n}d(t))].
\end{align}
Next, expressing the dynamics in the $\mathbf{z}$-coordinates and $\tau$-time, applying the time-scale transformation~\eqref{Eq:TimeScaleTransf}, and noticing that
\begin{equation}
    \frac{\mathrm{d}\mathbf{z}}{\mathrm{d}\tau}=\left.\frac{\mathrm{d}\mathbf{z}}{\mathrm{d}t}\right|_{t=\eta^{-1}T_c\left(1-\exp(-\alpha \tau)\right)}\frac{\mathrm{d}t}{\mathrm{d}\tau}
\end{equation}
where $\frac{\mathrm{d}t}{\mathrm{d}\tau}$ is given by~\eqref{Eq:dtdtau}, yields system~\eqref{Eq:AuxSystemAOS} with
$$
\pi(\tau):=\beta^{-1} (\alpha\eta^{-1}T_c)^{n-\rho}\exp(-\alpha(n-\rho)\tau)d(\varphi^{-1}(\tau)).
$$
Notice that $\pi(\tau)$ satisfies~\eqref{Eq:PiAux}. Since the origin of system~\eqref{Eq:AuxSystemAOS} is asymptotically stable and its settling-time function satisfies~\eqref{Eq:STFAOS}, then the settling-time of the closed-loop system under the controller~\eqref{Eq:POC} is given by~\eqref{Eq:Settling-timeAOS}, which completes the proof.
\qed
\end{proof}

\begin{corollary}
Let $\mathbf{v}(\mathbf{z})$ be an auxiliary controller such that the closed-loop system of~\eqref{Eq:AuxSystemAOS} is asymptotically stable with settling time function $\mathcal{T}(\mathbf{z}_0)$ given by~\eqref{Eq:STFAOS}. Then, under the controller~\eqref{Eq:POC}:
\begin{enumerate}
    \item if $\mathcal{T}(\mathbf{z}_0)=\infty$, then, the settling-time function~\eqref{Eq:Settling-timeAOS} satisfies $T(\mathbf{x}_0)=T_c$. Thus, controller~\eqref{Eq:POC} is a prescribed-time controller but $\lim_{t\to T(\mathbf{x}_0)^-}\kappa(t)=\infty$.
    \item if $\mathcal{T}(\mathbf{z}_0)<\infty$, but $T_f=\infty$, then, the time-varying gain $\kappa(T(\mathbf{x}_0))$ is finite but an unbounded function of the initial condition $\mathbf{x}_0$, i.e. $\sup_{\mathbf{x}_0\in\mathbb{R}^n}\kappa(T(\mathbf{x}_0))=\infty$. Thus, controller~\eqref{Eq:POC} is a predefined-time controller with a finite (but unbounded) time-varying gain at the settling time.
    \item if $\sup_{\mathbf{z}_0\in\mathbb{R}^n}\mathcal{T}(\mathbf{z}_0)\leq T_f$, with a known $T_f<\infty$, then, the settling-time function~\eqref{Eq:Settling-timeAOS} satisfies $T(\mathbf{x}_0)\leq T_c$, and $\kappa(t)$ is bounded by the $\kappa_{\max}:=\frac{\exp(\alpha T_f)-1}{\alpha T_c}$ for $t\in[0,T_c)$. Thus, controller~\eqref{Eq:POC} is a predefined-time controller with bounded time-varying gain.
    
\end{enumerate}
\end{corollary}

\begin{example}
Let $n=3$ and assume that for parameters $\alpha\geq 0$, $\rho=n$ and $T_f\in\bar{\mathbb{R}}_+$, an admissible auxiliary controller $v(\mathbf{z})=\omega(z_1,z_2,z_3)$ is given (fixed-time autonomous controllers with estimation of the settling time have been presented in~\cite{Zimenko2018,Cao2020PrespecifiableForm}). Thus, matrices $
\mathbf{Q}_\rho$ and $\mathbf{K}_{\rho}(t)$ are computed as
$$
\mathbf{Q}_\rho=
\begin{bmatrix}
1 & 0 & 0 \\
3\alpha & 1 &0 \\
9\alpha^2 & 5\alpha & 1
\end{bmatrix}
\text{ and }
\mathbf{K}_{\rho}(t)=\textrm{diag}(\kappa(t)^{-3},\kappa(t)^{-2},\kappa(t)^{-1}),
$$
respectively. Thus,
$$
\mathbf{Q}_{\rho} \mathbf{K}_{\rho}^{-1}(t) \mathbf{x}=
\begin{bmatrix}
\kappa(t)^{3}x_1 \\
3\alpha\kappa(t)^{3}x_1+\kappa(t)^{2}x_2\\
9\alpha^2\kappa(t)^{3}x_1+5\alpha\kappa(t)^{2}x_2+\kappa(t)x_3
\end{bmatrix}
$$
and
$$
\mathbf{b}_{1}^{T}(\mathbf{J}-\alpha \mathbf{D}_{\rho})^{n} \mathbf{K}_{\rho}^{-1}(t) \mathbf{x}=27\alpha^3\kappa(t)^{3}x_1+19\alpha^2\kappa(t)^{2}x_2+6\alpha\kappa(t)x_3.
$$
Therefore, we obtain the predefined-time controller~\eqref{Eq:POC} with $\phi(\mathbf{x},t;T_c)$ given by:
\begin{multline}
\phi(\mathbf{x},t;T_c)=\omega\left(\kappa(t)^{3}x_1,3\alpha\kappa(t)^{3}x_1+\kappa(t)^{2}x_2,9\alpha^2\kappa(t)^{3}x_1+5\alpha\kappa(t)^{2}x_2+\kappa(t)x_3\right)\\-27\alpha^3\kappa(t)^{3}x_1-19\alpha^2\kappa(t)^{2}x_2-6\alpha\kappa(t)x_3
\end{multline}
\end{example}

\subsection{Uniform Lyapunov stability of predefined-time controllers}

Since our approach for predefined-time control uses \emph{time-varying} gains to redesign an admissible auxiliary controller, it is essential to study the uniform (with respect to time) stability. This property has robustness implications, for instance, with respect to measurement noise, quantization, etc. We say that the origin of system~\eqref{Eq:IntegratorChain} is \emph{uniformly Lyapunov stable} \cite[Definition 4.4]{Khalil2002NonlinearSystems}, if for every $\epsilon > 0$ there exists $\delta > 0$ such that for all $s \ge 0$ and every admissible disturbance, $||\mathbf{x}(s)|| \le \delta$ implies $||\mathbf{x}(t)|| \le \epsilon$ for all $t \ge s$. 

\index{Measurement/Quatization noise}
\begin{example}
\label{Ex:LackofUniformity}
Let $\Delta=0$, $\rho=0$, $\alpha=1$, $T_f=\infty$ and  $v(\mathbf{z})=-18z_1-9z_2$. Notice that such $v(\mathbf{z})$ is an admissible auxiliary controller, since under such controller the state trajectory of the auxiliary system~\eqref{Eq:AuxSystemAOS} is:
\begin{align}
z_1(\tau;\mathbf{z}_0)&=\left(2 \exp(-3 t) - \exp(-6 t)\right)z_1(0)+\left(\frac{1}{3}\exp(-3 t) - \frac{1}{3} \exp(-6 t)\right)z_2(0)\\
z_2(\tau;\mathbf{z}_0)&=(6 \exp(-6 t) - 6 \exp(-3 t))z_1(0)+ (2 \exp(-6 t) - \exp(-3 t))z_2(0)
\end{align}
the origin of~\eqref{Eq:AuxSystemAOS} is asymptotically stable with $\mathcal{T}(\mathbf{z}_0)=\infty$, for every nonzero initial condition $\mathbf{z}_0$. Moreover, it is easy to verify that \eqref{Eq:Limit} holds. 

Thus, $$v(\beta^{-1}\mathbf{Q}_{\rho} \mathbf{K}_{\rho}^{-1}(t) \mathbf{x})=-18x_1-9\kappa(t)^{-1}x_2$$ and $$\beta^{-1}\mathbf{b}_{1}^{T}(\mathbf{J}-\alpha \mathbf{D}_{\rho})^{n} \mathbf{K}_{\rho}^{-1}(t) \mathbf{x}=-\kappa(t)^{-1}x_2$$ and using such admissible auxiliary controller we obtain, based on Theorem~\ref{Th:Main}, the prescribed-time controller:
\begin{equation}
\label{Eq:PrescribedSOC}
\phi(\mathbf{x},t;T_c)=-18\kappa(t)^2x_1-8\kappa(t)x_2.   
\end{equation}
The trajectory of such prescribed-time controller, for $T_c=10$ and initial condition $x_1(0)=x_2(0)=10$ is shown in the first column of Fig.~\ref{Fig:NonUniformity}. 

Next, consider a disturbance 
\begin{equation}
\label{Eq:PulseDist}
\mu(t;t_d)=\left\lbrace 
\begin{array}{lll}
     0.1& \text{if} & t_d\leq t<t_d+0.001  \\
     0 & & \text{otherwise} 
\end{array}
\right.
\end{equation}
where $t_d<T_c$, such that the prescribed-time controller becomes:
$$\phi\left(\mathbf{x}+
\begin{bmatrix}
\mu(t;t_d) \\ 0
\end{bmatrix},t;T_c\right),$$
with $\phi(\mathbf{x},t;T_c)$ as in~\eqref{Eq:PrescribedSOC}.

The second and third column of Fig.~\ref{Fig:NonUniformity} show the disturbance $\mu(t;t_d)$ (which could occur due to quantization, noise, etc) and the trajectory for $x_2(t)$, respectively, for different selections of $t_d$. As can be observed, an arbitrarily large transient can be obtained if $t_d$ is sufficiently close to $T_c$, which shows absence of uniform stability and a lack of robustness to disturbances in the state.

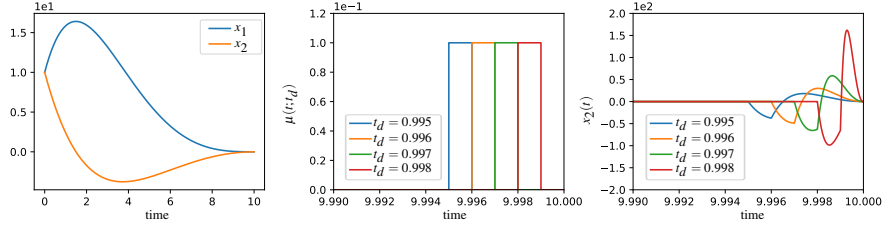
\begin{figure}
    \centering
\def\svgwidth{12cm}
\begingroup%
  \makeatletter%
  \providecommand\color[2][]{%
    \errmessage{(Inkscape) Color is used for the text in Inkscape, but the package 'color.sty' is not loaded}%
    \renewcommand\color[2][]{}%
  }%
  \providecommand\transparent[1]{%
    \errmessage{(Inkscape) Transparency is used (non-zero) for the text in Inkscape, but the package 'transparent.sty' is not loaded}%
    \renewcommand\transparent[1]{}%
  }%
  \providecommand\rotatebox[2]{#2}%
  \newcommand*\fsize{\dimexpr\f@size pt\relax}%
  \newcommand*\lineheight[1]{\fontsize{\fsize}{#1\fsize}\selectfont}%
  \ifx\svgwidth\undefined%
    \setlength{\unitlength}{855.59686279bp}%
    \ifx\svgscale\undefined%
      \relax%
    \else%
      \setlength{\unitlength}{\unitlength * \real{\svgscale}}%
    \fi%
  \else%
    \setlength{\unitlength}{\svgwidth}%
  \fi%
  \global\let\svgwidth\undefined%
  \global\let\svgscale\undefined%
  \makeatother%
  \begin{picture}(1,0.25993148)%
    \lineheight{1}%
    \setlength\tabcolsep{0pt}%
    \put(0,0){\includegraphics[width=\unitlength]{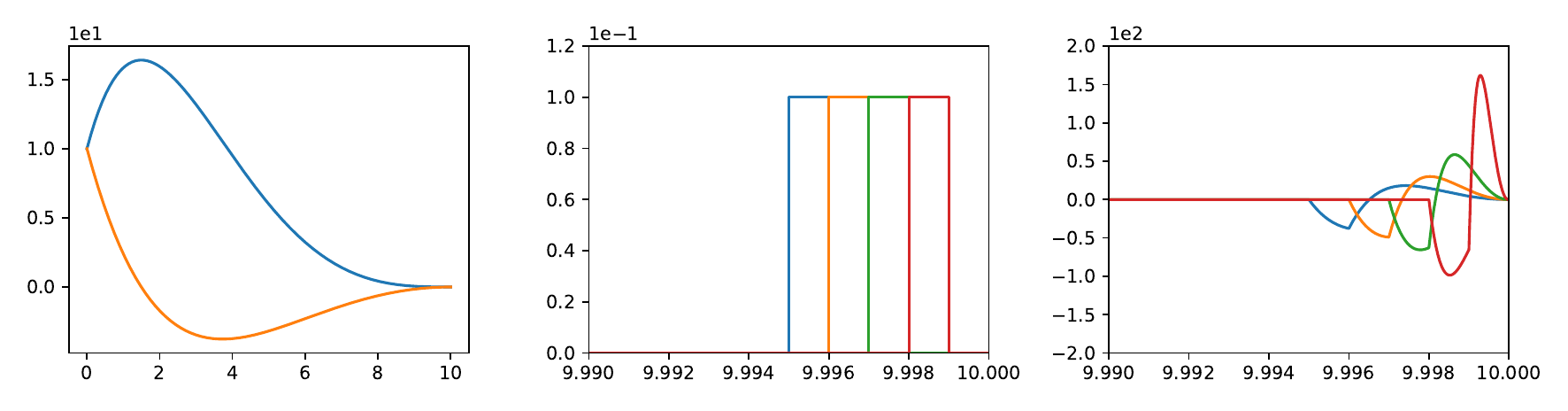}}%
    \put(0,0){\includegraphics[width=\unitlength]{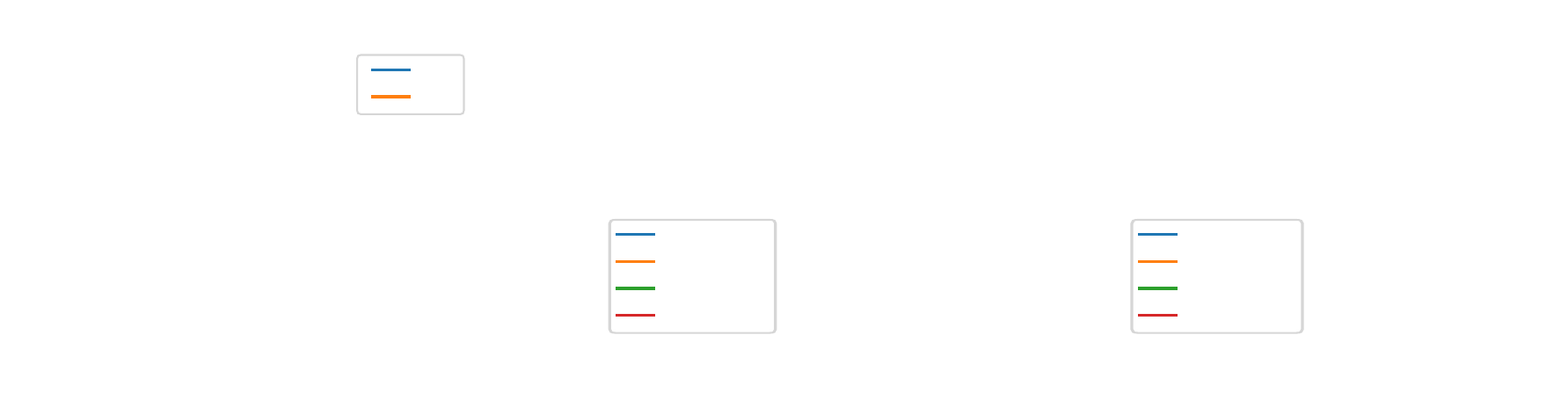}}%
    \tiny{
    \put(0.16715615,0.00383024){\makebox(0,0)[lt]{\lineheight{1.25}\smash{\begin{tabular}[t]{l}time\end{tabular}}}}%
    \put(0.82657884,0.00383024){\makebox(0,0)[lt]{\lineheight{1.25}\smash{\begin{tabular}[t]{l}time\end{tabular}}}}%
    \put(0.49686751,0.00383024){\makebox(0,0)[lt]{\lineheight{1.25}\smash{\begin{tabular}[t]{l}time\end{tabular}}}}%
    \put(0.2662081,0.2120877){\makebox(0,0)[lt]{\lineheight{1.25}\smash{\begin{tabular}[t]{l}$x_1$\end{tabular}}}}%
    \put(0.2662081,0.19455608){\makebox(0,0)[lt]{\lineheight{1.25}\smash{\begin{tabular}[t]{l}$x_2$\end{tabular}}}}%
    \put(0.42,0.10738305){\makebox(0,0)[lt]{\lineheight{1.25}\smash{\begin{tabular}[t]{l}$t_d=0.995$\end{tabular}}}}%
    \put(0.42,0.08985143){\makebox(0,0)[lt]{\lineheight{1.25}\smash{\begin{tabular}[t]{l}$t_d=0.996$\end{tabular}}}}%
    \put(0.42,0.07231981){\makebox(0,0)[lt]{\lineheight{1.25}\smash{\begin{tabular}[t]{l}$t_d=0.997$\end{tabular}}}}%
    \put(0.42,0.05654135){\makebox(0,0)[lt]{\lineheight{1.25}\smash{\begin{tabular}[t]{l}$t_d=0.998$\end{tabular}}}}%
    \put(0.755,0.10738305){\makebox(0,0)[lt]{\lineheight{1.25}\smash{\begin{tabular}[t]{l}$t_d=0.995$\end{tabular}}}}%
    \put(0.755,0.08985143){\makebox(0,0)[lt]{\lineheight{1.25}\smash{\begin{tabular}[t]{l}$t_d=0.996$\end{tabular}}}}%
    \put(0.755,0.07231981){\makebox(0,0)[lt]{\lineheight{1.25}\smash{\begin{tabular}[t]{l}$t_d=0.997$\end{tabular}}}}%
    \put(0.755,0.05654135){\makebox(0,0)[lt]{\lineheight{1.25}\smash{\begin{tabular}[t]{l}$t_d=0.998$\end{tabular}}}}%
    \put(0.33,0.10913306){\rotatebox{90}{\makebox(0,0)[lt]{\lineheight{1.25}\smash{\begin{tabular}[t]{l}$\mu(t;t_d)$\end{tabular}}}}}%
    \put(0.65925908,0.10913306){\rotatebox{90}{\makebox(0,0)[lt]{\lineheight{1.25}\smash{\begin{tabular}[t]{l}$x_2(t)$\end{tabular}}}}}%
    }
  \end{picture}%
\endgroup%
    \caption{Simulation of Example~\ref{Ex:LackofUniformity}, showing the lack of robustness to measurement noise of a prescribed-time algorithm. On the left, the behavior of a prescribed control with $T_c=10$ and without disturbance. In the center, a set of pulse disturbances in~\eqref{Eq:PulseDist}. On the right, the behavior of the closed-loop system under the prescribed control and in the presence of disturbance~\eqref{Eq:PulseDist}.}
    \label{Fig:NonUniformity}
\end{figure}
\end{example}

Note that when $n=1$, the change of variables between $\mathbf{z}$ and $\mathbf{x}$ in \eqref{eq:transformation} does not depend on time. Thus, if the dynamics of $\mathbf{z}$ in \eqref{Eq:AuxSystem} is uniformly Lyapunov stable, this property is transferred to $\mathbf{x}$ by continuity, as well. However, this reasoning fails with $n>1$ due to the time dependence of \eqref{eq:transformation}. In the following, we establish under which conditions, uniform Lyapunov stability is attained for the origin of the closed-loop system~\eqref{Eq:IntegratorChain} with the controller in \eqref{Eq:POC}. 

\index{Uniform Lyapunov stability}

\begin{proposition}
\label{Prop:Uniformity}
Consider $n>1$ and assume that the origin of the system \eqref{Eq:AuxSystem} is uniformly Lyapunov stable and $0\leq\rho<n-1$. Then, the origin of the closed-loop system~\eqref{Eq:IntegratorChain} with the controller in \eqref{Eq:POC} is uniformly Lyapunov stable if and only if $\kappa(t)$ is uniformly bounded. 
\end{proposition}
\begin{proof}
First, note that uniform Lyapunov stability for $t\geq T_c$ of \eqref{Eq:IntegratorChain} follows since the controller \eqref{Eq:POC} is time-independent for such $t$, and $w(\mathbf{x}, \Delta)$ is assumed to make the origin of the system stable. Thus, we only need to analyse the uniform Lyapunov stability property of \eqref{Eq:IntegratorChain} for $t\in[0,T_c)$. Note that if the solution $\mathbf{z}(\tau)$ of \eqref{Eq:AuxSystem} is uniformly Lyapunov stable, then the same applies to $\mathbf{z}(\varphi(t))$, $\forall t\in[0,T_c)$. Recall that the solution $\mathbf{x}(t)$ of \eqref{Eq:IntegratorChain} and $\mathbf{z}(\varphi(t))$ are related through the transformation in \eqref{eq:transformation} for $t\in[0,T_c)$. In addition, note that due to stability property of trajectories $\mathbf{x}(t)$ obtained from Theorem \ref{Th:Main}, i.e. that $\mathbf{x}(T_c)=0$ regardless of the initial conditions and disturbances, we can continue trajectories from $t=T_c$.

We start by showing that $\kappa(t)$ uniformly bounded implies uniform Lyapunov stability of the origin of system~\eqref{Eq:IntegratorChain}. Let $0\leq s< t<T_c$ and use Rayleigh's inequality in \eqref{eq:transformation} to obtain:
\begin{equation}
\label{eq:rayleigh}
\beta\underline{\sigma}(\mathbf{Q}_{\rho}^{-1})\underline{\sigma}(\Lambda_\rho(t))\|\mathbf{z}(\varphi(t))\|\leq \|\mathbf{x}(t)\|\leq\beta\overline{\sigma}(\mathbf{Q}_{\rho}^{-1})\overline{\sigma}(\Lambda_\rho(t))\|\mathbf{z}(\varphi(t))\|
\end{equation}
where $\underline{\sigma}(\bullet), \overline{\sigma}(\bullet)$ denote minimum and maximum singular values respectively. In addition, note that uniformly bounded $\kappa(t)$ implies that there exists  $0<\underline{\lambda},\overline{\lambda}\in\mathbb{R}$ such that 
\begin{equation}
\begin{array}{lcr}
\underline{\lambda}<&\underline{\sigma}(\mathbf{K}_{\rho}(t))=\min\{\kappa(t)^{-\rho},\dots,\kappa(t)^{n-\rho-1}\}& \\ &\overline{\sigma}(\mathbf{K}_{\rho}(t))=\max\{\kappa(t)^{-\rho},\dots,\kappa(t)^{n-\rho-1}\}&<\overline{\lambda}, 
\end{array} 
\end{equation}
for any $t\in[0,T_c]$. Thus, \eqref{eq:rayleigh} becomes:

\begin{equation}
\label{eq:rayleigh2}
\beta\underline{\sigma}(\mathbf{Q}_{\rho}^{-1})\underline{\lambda}\|\mathbf{z}(\varphi(t))\|\leq \|\mathbf{x}(t)\|\leq\beta\overline{\sigma}(\mathbf{Q}_{\rho}^{-1})\overline{\lambda}\|\mathbf{z}(\varphi(t))\|.
\end{equation}

Now, choose any $\epsilon>0$ and let $\epsilon_z = \epsilon\left(\beta\overline{\sigma}(\mathbf{Q}_\rho^{-1})\overline{\lambda}\right)^{-1}$ and note that $\epsilon_z>0$ since $\overline{\lambda}<+\infty$. For such $\epsilon_z>0$, there exists $\delta_z>0$ such that $\|\mathbf{z}(\varphi(s))\|\leq \delta_z$ implies $\|\mathbf{z}(\varphi(t))\|\leq \epsilon_z$ for $\varphi(t)\geq \varphi(s)$ due to Lyapunov stablility and time invariance (and hence, uniform Lyapunov stability) of \eqref{Eq:AuxSystem}. Thus, in order to see how this property is transferred to $\mathbf{x}$, let $\delta = \delta_z\beta\underline{\sigma}(\mathbf{Q}_\rho^{-1})\underline{\lambda}$. Hence, using the first inequality in \eqref{eq:rayleigh2} it is obtained that $\|\mathbf{x}(s)\|\leq \delta= \delta_z\beta\underline{\sigma}(\mathbf{Q}_\rho^{-1})\underline{\lambda}$ implies $\|\mathbf{z}(\varphi(s))\|\leq \delta_z$. This in turn implies $\|\mathbf{z}(\varphi(t))\|\leq \epsilon_z$.
Using the second inequality in \eqref{eq:rayleigh2} we obtain  $\|\mathbf{x}(t)\|\leq \beta\overline{\sigma}(\mathbf{Q}_\rho^{-1})\overline{\lambda} \epsilon_{z} = \epsilon$, showing uniform Lyapunov stability for the origin of the closed-loop system~\eqref{Eq:IntegratorChain} with the controller in \eqref{Eq:POC}.

Next, we show that if $\kappa(t)$ is not bounded, then the origin of the closed-loop system~\eqref{Eq:IntegratorChain} is not uniformly Lyapunov stable for $\rho \in[0, n-1)$. In particular, we will show that for any $\delta,\epsilon>0$, there exist $s,t$ with $0\leq s<t\leq T_c$, an admissible disturbance and a trajectory $x$ of \eqref{Eq:IntegratorChain} which satisfies both $\|\mathbf{x}(s)\|\leq \delta$ and $\|\mathbf{x}(t)\|> \epsilon$, which is the direct negation of the uniform Lyapunov stability statement. This means that we only need to find a single trajectory of \eqref{Eq:IntegratorChain} which fails to fulfill the uniformity bounds. 

In this sense, we can focus on $\Delta=0$ and consider, for fixed $\delta$ and arbitrary $\tau_0\ge 0$, any trajectory $\mathbf{z}(\tau)$ of \eqref{Eq:AuxSystem} with $\mathbf{z}(\tau_0) = w\mathbf{Q}_\rho \mathbf{b}_{1}$ with non-zero constant $w$ with $|w|\leq \delta/(\beta \kappa(0)^{-\rho})$. The proof strategy is to show that one of the zero components of $\mathbf{x}(s)$ with $\tau_0=\varphi(s)$ need to increase in magnitude at some later time $s<t$, enough to make $\mathbf{x}(t)$ as large as desired thanks to the unboundedness of the gain.

To guarantee that one component of $\mathbf{x}(s)$ cannot remain at zero, we show that there is no vector $\mathbf{h}=[h_1,\dots,h_n]^T\in\mathbb{R}^n$ with $h_n=0$ such that $\left.\frac{\mathrm{d}}{\mathrm{d}\tau}\mathbf{Q}_\rho^{-1} \mathbf{z}(\tau)\right|_{\tau=\tau_0}= \mathbf{h}$ for this trajectory, meaning that the last component of $\mathbf{x}(s)$ cannot remain at the origin. Assume the opposite, which implies that
$$
\mathbf{h}=\left.\mathbf{Q}_\rho^{-1}\frac{\mathrm{d}z}{\mathrm{d}\tau}\right|_{\tau=\tau_0} =  \mathbf{Q}_\rho^{-1}\mathbf{b}_{n}{v}(\mathbf{z}(\tau_0)) + w\mathbf{Q}_\rho^{-1}\mathbf{J}\mathbf{Q}_\rho \mathbf{b}_{1}=\mathbf{b}_{n}{v}(\mathbf{z}(\tau_0))
$$
since $\mathbf{Q}_\rho^{-1}\mathbf{b}_{n}=\mathbf{b}_{n}$ and $\mathbf{Q}_\rho^{-1}\mathbf{J}\mathbf{Q}_\rho \mathbf{b}_{1}=0$. However, this is impossible since $h_n=0$ but ${v}(\mathbf{z}(\tau_0))\neq 0$. Therefore,  $\left.\frac{\mathrm{d}}{\mathrm{d}\tau}\mathbf{b}_{n}^T\mathbf{Q}_\rho^{-1}\mathbf{z}(\tau)\right|_{\tau=\tau_0}$ is non zero. The previous argument, in addition to the fact that \eqref{Eq:AuxSystem} is time-invariant and ${v}(\bullet)$ is continuous at $\mathbf{z}(\tau_0)$, implies that there exist positive constants $\tilde{\tau}$, $\tilde{\epsilon}$, which only depend on $\delta$, such that
$|\mathbf{b}_{n}^T \mathbf{Q}_\rho^{-1} \mathbf{z}(\tau_0 + \tilde{\tau})|> \tilde{\epsilon}$. 

Select now $s\geq 0$ such that $\beta \kappa(s)^{n-\rho-1} \tilde{\epsilon} > \epsilon$ which is possible since $\kappa(\bullet)$ is unbounded and $n-\rho-1>0$. Set $\tau_0 = \varphi(s)$ and note that from \eqref{eq:transformation}, one then obtains
\begin{align}
    \mathbf{x}(s) &= \beta \mathbf{K}_{\rho}(s) \mathbf{Q}_\rho^{-1} \mathbf{z}(\varphi(s))= \beta \mathbf{K}_{\rho}(s) \mathbf{Q}_\rho^{-1}\mathbf{z}(\tau_0)\\ &= \beta w \mathbf{K}_{\rho}(s) \mathbf{b}_{1} = \beta w \kappa(s)^{-\rho} \mathbf{b}_{1}
\end{align}
and consequently $\|\mathbf{x}(s)\|= \beta |w|\kappa(s)^{-\rho}\leq \beta |w|\kappa(0)^{-\rho}  \le \delta$. Moreover, one has for $t = \varphi^{-1}(\varphi(s)+\tilde{\tau}) <T_c$
\begin{align}
    x_n(t) &= \mathbf{b}_{n}^T \mathbf{x}(t)= \beta \kappa(t)^{n-\rho-1} \mathbf{b}_{n}^T \mathbf{Q}_\rho^{-1} \mathbf{z}(\varphi(t)) \\ & = \beta \kappa(t)^{n-\rho-1} \mathbf{b}_{n}^T \mathbf{Q}_\rho^{-1} \mathbf{z}(\tau_0 + \tilde{\tau})
\end{align}
and hence $\|\mathbf{x}(t)\| \ge |x_n(t)| \ge \beta \kappa(t)^{n-\rho-1} \tilde{\epsilon} \ge \beta \kappa(s)^{n-\rho-1} \tilde{\epsilon} > \epsilon$, since $\kappa(t)$ is non decreasing and $n-\rho-1>0$, completing the proof.
\qed
\end{proof}

Proposition~\ref{Prop:Uniformity} implies that, with the proposed approach, no uniformly stable \emph{prescribed-time} controller can be obtained for a chain of integrators of order greater than one. Ensuring uniform stability for \emph{predefined-time} controllers, on the other hand, may be achieved by ensuring that the time-varying gain stays bounded. Although this proof is particular for our approach, the proof suggests that we can take a similar path to show the non-uniformity of other prescribed-time control methods. Since, as illustrated in Example~\ref{Ex:LackofUniformity}, the lack of uniform stability implies an inherent lack of robustness, then studying the uniform stability property in the prescribed-time control literature is essential.

\index{Measurement/Quatization noise}
\begin{example}
\label{Ex:Uniformity}
Let us revisit the controller in Example~\ref{Ex:SOSAut} for a perturbed system with disturbance $d(t)=\sin(t)$. Similarly as in Example~\ref{Ex:LackofUniformity} consider the disturbance $\mu(t;t_d)$ in~\eqref{Eq:PulseDist}, such that the predefined-time controller becomes $$\phi\left(\mathbf{x}+
\begin{bmatrix}
\mu(t;t_d) \\ 0
\end{bmatrix},t;T_c\right).$$
The second and third columns of Fig.~\ref{Fig:Uniformity} show the disturbance $\mu(t;t_d)$ (which could occur due to quantization, noise, etc) and the trajectory for $x_2(t)$, respectively, for different selections of $t_d$. As can be observed, contrary to the case in Example~\ref{Ex:LackofUniformity}, the transient behavior due to the perturbation remains bounded, no matter how close to $T_c$ the disturbance $\mu(t;t_d)$ occurs. Notice that according to Proposition~\ref{Prop:Uniformity} this predefined-time controller is uniformly Lyapunov stable since the time-varying gain is bounded.

\begin{figure}
    \centering
\def\svgwidth{12cm}
\begingroup%
  \makeatletter%
  \providecommand\color[2][]{%
    \errmessage{(Inkscape) Color is used for the text in Inkscape, but the package 'color.sty' is not loaded}%
    \renewcommand\color[2][]{}%
  }%
  \providecommand\transparent[1]{%
    \errmessage{(Inkscape) Transparency is used (non-zero) for the text in Inkscape, but the package 'transparent.sty' is not loaded}%
    \renewcommand\transparent[1]{}%
  }%
  \providecommand\rotatebox[2]{#2}%
  \newcommand*\fsize{\dimexpr\f@size pt\relax}%
  \newcommand*\lineheight[1]{\fontsize{\fsize}{#1\fsize}\selectfont}%
  \ifx\svgwidth\undefined%
    \setlength{\unitlength}{855.59686279bp}%
    \ifx\svgscale\undefined%
      \relax%
    \else%
      \setlength{\unitlength}{\unitlength * \real{\svgscale}}%
    \fi%
  \else%
    \setlength{\unitlength}{\svgwidth}%
  \fi%
  \global\let\svgwidth\undefined%
  \global\let\svgscale\undefined%
  \makeatother%
  \begin{picture}(1,0.25993148)%
    \lineheight{1}%
    \setlength\tabcolsep{0pt}%
    \put(0,0){\includegraphics[width=\unitlength]{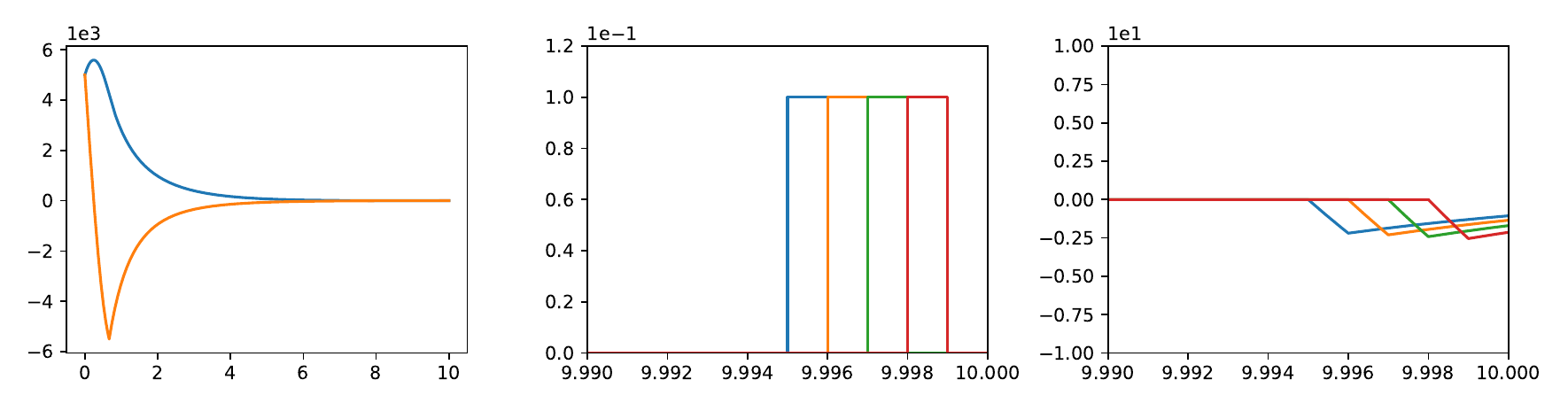}}%
    \put(0,0){\includegraphics[width=\unitlength]{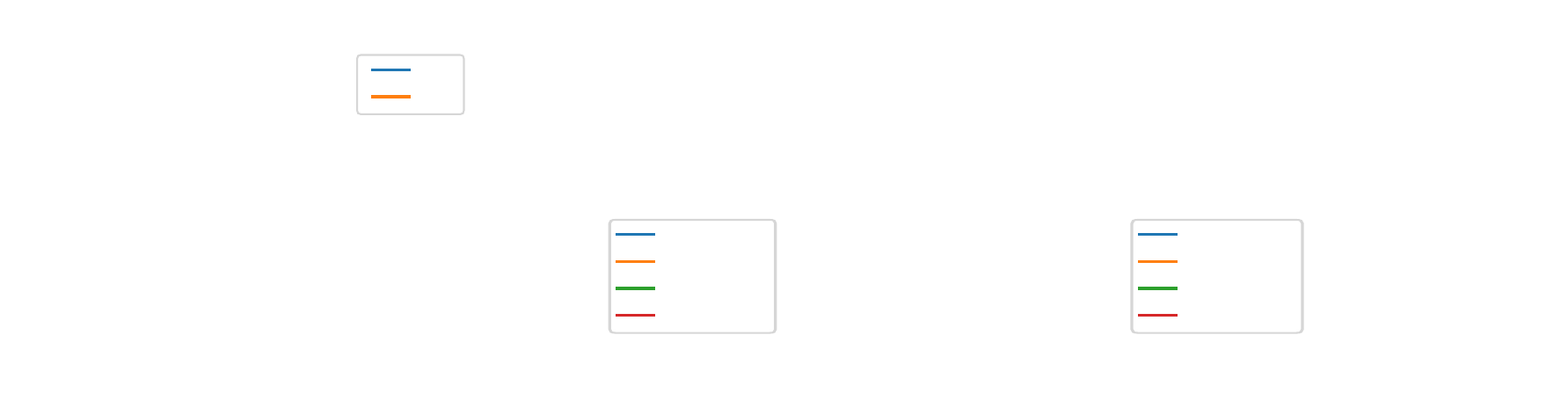}}%
    \tiny{
    \put(0.16715615,0.00383024){\makebox(0,0)[lt]{\lineheight{1.25}\smash{\begin{tabular}[t]{l}time\end{tabular}}}}%
    \put(0.82657884,0.00383024){\makebox(0,0)[lt]{\lineheight{1.25}\smash{\begin{tabular}[t]{l}time\end{tabular}}}}%
    \put(0.49686751,0.00383024){\makebox(0,0)[lt]{\lineheight{1.25}\smash{\begin{tabular}[t]{l}time\end{tabular}}}}%
    \put(0.2662081,0.2120877){\makebox(0,0)[lt]{\lineheight{1.25}\smash{\begin{tabular}[t]{l}$x_1$\end{tabular}}}}%
    \put(0.2662081,0.19455608){\makebox(0,0)[lt]{\lineheight{1.25}\smash{\begin{tabular}[t]{l}$x_2$\end{tabular}}}}%
    \put(0.42,0.10738305){\makebox(0,0)[lt]{\lineheight{1.25}\smash{\begin{tabular}[t]{l}$t_d=0.995$\end{tabular}}}}%
    \put(0.42,0.08985143){\makebox(0,0)[lt]{\lineheight{1.25}\smash{\begin{tabular}[t]{l}$t_d=0.996$\end{tabular}}}}%
    \put(0.42,0.07231981){\makebox(0,0)[lt]{\lineheight{1.25}\smash{\begin{tabular}[t]{l}$t_d=0.997$\end{tabular}}}}%
    \put(0.42,0.05654135){\makebox(0,0)[lt]{\lineheight{1.25}\smash{\begin{tabular}[t]{l}$t_d=0.998$\end{tabular}}}}%
    \put(0.755,0.10738305){\makebox(0,0)[lt]{\lineheight{1.25}\smash{\begin{tabular}[t]{l}$t_d=0.995$\end{tabular}}}}%
    \put(0.755,0.08985143){\makebox(0,0)[lt]{\lineheight{1.25}\smash{\begin{tabular}[t]{l}$t_d=0.996$\end{tabular}}}}%
    \put(0.755,0.07231981){\makebox(0,0)[lt]{\lineheight{1.25}\smash{\begin{tabular}[t]{l}$t_d=0.997$\end{tabular}}}}%
    \put(0.755,0.05654135){\makebox(0,0)[lt]{\lineheight{1.25}\smash{\begin{tabular}[t]{l}$t_d=0.998$\end{tabular}}}}%
    \put(0.33,0.10913306){\rotatebox{90}{\makebox(0,0)[lt]{\lineheight{1.25}\smash{\begin{tabular}[t]{l}$\mu(t;t_d)$\end{tabular}}}}}%
    \put(0.65925908,0.10913306){\rotatebox{90}{\makebox(0,0)[lt]{\lineheight{1.25}\smash{\begin{tabular}[t]{l}$x_2(t)$\end{tabular}}}}}%
    }
  \end{picture}%
\endgroup%
    \caption{Simulation of Example~\ref{Ex:Uniformity}, showing robustness to measurement noise of a prescribed-time algorithm with bounded time varying gains. On the left, the behavior of a prescribed control with $T_c=10$ and without disturbance. In the center, a set of pulse disturbances in~\eqref{Eq:PulseDist}. On the right, the behavior of the closed-loop system under the prescribed control and in the presence of disturbance~\eqref{Eq:PulseDist}.}
    \label{Fig:Uniformity}
\end{figure}
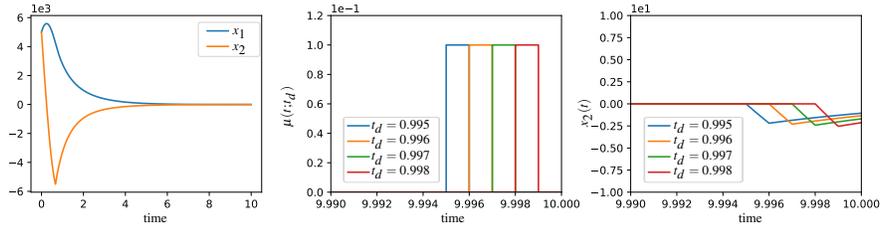

\end{example}

\section{Conclusion}

This chapter presents a methodology to design robust controllers achieving fixed-time stability with a desired upper bound for the settling time (\textit{UBST}). We show that the closed-loop system under the proposed controller is related to a suitable auxiliary system through a time-varying coordinate change and a time-scale transformation. The methodology is motivated by the analysis of the first-order case, where application to a linear controller leads to a minimum energy solution and generally allows to reduce control energy when redesigning other controllers.
Depending on the convergence properties of the resulting auxiliary system, interesting features are obtained in the closed-loop system. For instance, obtaining a prescribed-time controller steers the state to the origin in the desired time, regardless of the initial condition, but with time-varying gains that tend to infinity. Alternatively, we present conditions under which a predefined-time controller is obtained with bounded time-varying gains.

Since the proposed controller is time-varying, it is essential to study its uniform stability properties. For this purpose, we show that uniform boundedness of the time-varying gain is necessary and sufficient for uniform Lyapunov stability of the closed-loop system obtained with our approach. It is moreover shown that such boundedness of the gain can be achieved by redesigning an existing fixed-time controller.

\begin{acknowledgement}
Work partially supported by the Christian Doppler Research Association, the Austrian Federal Ministry of Labour and Economy, and the National Foundation for Research, Technology and Development, by Agencia I+D+i grant PICT 2018-01385, Argentina and by Consejo Nacional de Ciencia y Tecnología (CONACYT-Mexico) scholarship with grant 739841.
\end{acknowledgement}
\section{Appendix}
\addcontentsline{toc}{section}{Appendix}

\subsection{Auxiliary lemmas}

Let us introduce the following Lemmas, on some properties of matrix $\mathbf{Q}_{\rho}$ and the time-varying matrix $\mathbf{K}_{\rho}(t)$. 

\begin{lemma}
\label{Lem:DQ:identities} 
Let $\mathbf{D}_{\rho}\in\mathbb{R}^{n\times n}$ and $\mathbf{Q}_{\rho}\in\mathbb{R}^{n\times n}$ be defined as in~\eqref{Eq:MatrixQp}. Then, $\mathbf{Q}_{\rho}\in\mathbb{R}^{n\times n}$ is a lower triangular matrix satisfying
\begin{equation}
 \label{Eq:QpId1}
\mathbf{J}+\mathbf{A} =\mathbf{Q}_{\rho}(\mathbf{J}-\alpha \mathbf{D}_{\rho})\mathbf{Q}_{\rho}^{-1}
\end{equation}
where 
\begin{equation}
\label{Eq:A}
    \mathbf{A}=\mathbf{b}_{n}\mathbf{b}_{1}^{T}(\mathbf{J}-\alpha \mathbf{D}_{\rho})^{n}\mathbf{Q}_{\rho}^{-1}
\end{equation} with $\mathbf{b}_{n}=[0,\cdots,0,1]^T\in\mathbb{R}^{n\times 1}$. 
\end{lemma}
\begin{proof}
Notice that by construction $\mathbf{Q}_{\rho}$ is a lower triangular matrix with ones over the diagonal. Moreover, $\mathbf{J}$ is an upper shift matrix, thus
\begin{equation}
\mathbf{J} \mathbf{Q}_{\rho}= \begin{bmatrix}
\mathbf{b}_{1}^{T}(\mathbf{J}-\alpha \mathbf{D}_{\rho})\\
\vdots \\ 
\mathbf{b}_{1}^{T}(\mathbf{J}-\alpha \mathbf{D}_{\rho})^{n-1}\\
 \mathbf{0}_n^T
\end{bmatrix} \text{ and } 
\mathbf{J} \mathbf{Q}_{\rho} +\mathbf{A} \mathbf{Q}_{\rho}= \begin{bmatrix}
\mathbf{b}_{1}^{T}(\mathbf{J}-\alpha \mathbf{D}_{\rho})\\
\vdots \\ 
\mathbf{b}_{1}^{T}(\mathbf{J}-\alpha \mathbf{D}_{\rho})^{n-1}\\
\mathbf{b}_{1}^{T}(\mathbf{J}-\alpha \mathbf{D}_{\rho})^{n}
\end{bmatrix}
\end{equation}
where $\mathbf{0}_n\in\mathbb{R}^n$ is a zero vector. 
Therefore, $
\mathbf{J} \mathbf{Q}_{\rho} -\mathbf{A} \mathbf{Q}_{\rho}=\mathbf{Q}_{\rho}(\mathbf{J}-\alpha \mathbf{D}_{\rho})
$ which completes the proof.
\qed
\end{proof}

\begin{lemma}
\label{Lem:kappa:identities}
Let $\kappa(t)$ be given as in~\eqref{Eq:kappa}, with $\eta$ as in~\eqref{Eq:Eta}, and let 
\begin{equation}
\mathbf{K}_{\rho}(t):=\textrm{diag}(\kappa(t)^{-\rho},\kappa(t)^{1-\rho},\ldots,\kappa(t)^{n-\rho-1}),
\end{equation}
where $\rho\in[0,n]$. Then, the following identities hold:
\begin{align}
   \frac{\mathrm{d}}{\mathrm{d}t}\mathbf{K}_{\rho}(t)^{-1} &= -\alpha \kappa(t)\mathbf{D}_{\rho}\mathbf{K}_{\rho}(t)^{-1} \label{Eq:kID1}\\
   \mathbf{K}_{\rho}^{-1}(t)\mathbf{J}\mathbf{K}_{\rho}(t)&=\kappa(t)\mathbf{J}.\label{Eq:kID2}
\end{align}
\end{lemma}
\begin{proof}

A direct calculation yields
\begin{align}
   \frac{\mathrm{d}}{\mathrm{d}t}\mathbf{K}_{\rho}(t)^{-1} &= \frac{\mathrm{d}}{\mathrm{d}t}\textrm{diag}(\kappa(t)^{\rho},\kappa(t)^{\rho-1},\ldots,\kappa(t)^{\rho-n+1})\\
   &=\dot{\kappa}(t)\kappa(t)^{-1}\textrm{diag}(\rho\kappa(t)^{\rho},\rho-1\kappa(t)^{\rho-1},\ldots,\rho-n+1\kappa(t)^{\rho-n+1}).
\end{align}
Since $\dot{\kappa}(t)\kappa(t)^{-1}=\alpha \kappa(t)$, equation~\eqref{Eq:kID1} follows trivially by definition of $\mathbf{D}_{\rho}$.

Now, to show that~\eqref{Eq:kID2} holds, notice that since $\mathbf{J}$ is an upper shift matrix. Thus,
\begin{align}
\mathbf{K}_{\rho}^{-1}(t)\mathbf{J}\mathbf{K}_{\rho}(t)&=\kappa(t)
\mathbf{K}_{\rho}^{-1}(t)
\begin{bmatrix}
0 & \kappa(t)^{-\rho} & 0 & \cdots &0 & 0\\
0 & 0 & \kappa(t)^{1-\rho} & \cdots &0 & 0\\
\vdots & \vdots & \vdots & \ddots & \vdots & \vdots\\
0 & 0 & 0 & \cdots & \kappa(t)^{n-3+\rho} & 0\\
0 & 0 & 0 & \cdots & 0 & \kappa(t)^{n-2+\rho}\\
0 & 0 & 0 & \cdots & 0 & 0\\
\end{bmatrix}\\
&=\kappa(t)\mathbf{J},
\end{align}
which completes the proof.
\qed
\end{proof}

\subsection{Some admissible auxiliary controllers}
\label{AppendixNonAut}

\begin{theorem}(\cite[Theorem~3]{Aldana-Lopez2018})
\label{th:tf_poly}
Consider a controller
\begin{equation}
\label{Eq:HomFixedPoly}
u = -\left[(a|x|^{p} + b|x|^{q})^k+\zeta\right] \sign{x},
\end{equation}
where $\zeta\geq\Delta$, $a,b,p,q,k>0$ are system parameters which satisfy the constraints $kp<1$, and $kq>1$. Then, the origin of~\eqref{Eq:FOS} under the controller~\eqref{Eq:HomFixedPoly} is fixed-time stable and the settling-time function satisfies $\sup_{x_0 \in \mathbb{R}} T(x_0)=\gamma$, where
\begin{equation}
\label{Eq:MinUpperEstimate}
\gamma=\frac{\Gamma \left(m_p\right) \Gamma \left(m_q\right)}{a^{k}\Gamma (k) (q-p)}\left(\frac{a}{b}\right)^{m_p},    
\end{equation}
with $m_p=\frac{1-k p}{q-p}$ and $m_q=\frac{k q-1}{q-p}$. 
\end{theorem}

\begin{theorem} (\cite[Theorem~4]{Aldana-Lopez2018})
\label{thm:socont} Consider a second-order perturbed chain of integrators, and let $a_1,a_2,b_1,b_2,p,q,k>0$, $kp<1$, $kq>1$, $T_{f_1},T_{f_2}>0$, $\zeta\geq\Delta$, and \[\gamma_1=\frac{\Gamma \left(\frac{1}{4}\right)^2 }{2a_1^{1/2}\Gamma\left(\frac{1}{2}\right)}\left(\frac{a_1}{b_1}\right)^{1/4},\text{ and } \gamma_2=\frac{\Gamma \left(m_{p}\right) \Gamma \left(m_{q}\right)}{a_2^{k}\Gamma (k) (q-p)}\left(\frac{a_2}{b_2}\right)^{m_{p}},\] with $m_{p}=\frac{1-kp}{q-p}$ and $m_{q}=\frac{kq-1}{q-p}$. If the control input is selected as
\begin{equation}\label{eq:uso}
u=-\left[\frac{\gamma_2}{T_{f_2}}\left(a_2\abs{\sigma}^{p}+b_2\abs{\sigma}^{q}\right)^{k}+\frac{\gamma_1^2}{2T_{f_1}^2}\left(a_1+3b_1x_1^2\right)+\zeta\right]\sign{\sigma},
\end{equation}
where the sliding variable $\sigma$ is defined as
\begin{equation}\label{eq:sigmaso}
\sigma=x_2+\barpow{\barpow{x_2}^2+\frac{2\gamma_1^2}{T_{f_1}^2}\left(a_1\barpow{x_1}^1+b_1\barpow{x_1}^3\right)}^{1/2},
\end{equation}
then the origin $(x_1,x_2)=(0,0)$ of system~\eqref{Eq:IntegratorChain}, with $n=2$, is fixed-time stable with \textit{UBST} given by $T_f=T_{f_1}+T_{f_2}$.
\end{theorem}

\bibliographystyle{spmpsci}

\begin{thebibliography}{10}
\providecommand{\url}[1]{{#1}}
\providecommand{\urlprefix}{URL }
\expandafter\ifx\csname urlstyle\endcsname\relax
  \providecommand{\doi}[1]{DOI~\discretionary{}{}{}#1}\else
  \providecommand{\doi}{DOI~\discretionary{}{}{}\begingroup
  \urlstyle{rm}\Url}\fi

\bibitem{Aldana-Lopez2018}
Aldana-L{\'{o}}pez, R., G{\'{o}}mez-Guti{\'{e}}rrez, D.,
  Jim{\'{e}}nez-Rodr{\'{i}}guez, E., S{\'{a}}nchez-Torres, J.D., Defoort, M.:
  {Enhancing the settling time estimation of a class of fixed-time stable
  systems}.
\newblock International Journal of Robust and Nonlinear Control
  \textbf{29}(12), 4135--4148 (2019).
\newblock \doi{10.1002/rnc.4600}

\bibitem{aldana2019design}
Aldana-L{\'{o}}pez, R., G{\'{o}}mez-Guti{\'{e}}rrez, D.,
  Jim{\'{e}}nez-Rodr{\'{i}}guez, E., S{\'{a}}nchez-Torres, J.D., Defoort, M.:
  {Generating new classes of fixed-time stable systems with predefined upper
  bound for the settling time}.
\newblock International Journal of Control \textbf{95}(10), 2802--2814 (2022).
\newblock \doi{10.1080/00207179.2021.1936190}.
\newblock \urlprefix\url{https://doi.org/10.1080/00207179.2021.1936190}

\bibitem{Andrieu2008}
Andrieu, V., Praly, L., Astolfi, A.: {Homogeneous approximation, recursive
  observer design, and output feedback}.
\newblock SIAM Journal on Control and Optimization \textbf{47}(4), 1814--1850
  (2008).
\newblock \doi{10.1137/060675861}

\bibitem{Andrieu2009HomogeneityDesign}
Andrieu, V., Praly, L., Astolfi, A.: {Homogeneity in the bi-limit as a tool for
  observer and feedback design}.
\newblock In: Proceedings of the IEEE Conference on Decision and Control, pp.
  1050--1055 (2009).
\newblock \doi{10.1109/CDC.2009.5400263}

\bibitem{Cao2020PrespecifiableForm}
Cao, Y., Wen, C., Tan, S., Song, Y.: {Prespecifiable fixed-time control for a
  class of uncertain nonlinear systems in strict-feedback form}.
\newblock International Journal of Robust and Nonlinear Control \textbf{30}(3),
  1203--1222 (2020).
\newblock \doi{10.1002/RNC.4820}

\bibitem{Chitour2020StabilizationTime}
Chitour, Y., Ushirobira, R., Bouhemou, H.: {Stabilization for a perturbed chain
  of integrators in prescribed time}.
\newblock SIAM Journal on Control and Optimization \textbf{58}(2), 1022--1048
  (2020).
\newblock \doi{10.1137/19M1285937}

\bibitem{Cruz-Zavala2021High-orderBi-limit}
Cruz-Zavala, E., Moreno, J.A.: {High-order sliding-mode control design
  homogeneous in the bi-limit}.
\newblock International Journal of Robust and Nonlinear Control \textbf{31}(9),
  3380--3416 (2021).
\newblock \doi{10.1002/RNC.5242}

\bibitem{Ding2016SimpleController}
Ding, S., Levant, A., Li, S.: {Simple homogeneous sliding-mode controller}.
\newblock Automatica \textbf{67}(5), 22--32 (2016).
\newblock \doi{10.1016/j.automatica.2016.01.017}

\bibitem{Filippov1988DifferentialSides}
Filippov, A.F.: {Differential equations with discontinuous righthand sides}.
\newblock Dordrecht: Kluwer Academic Publishers (1988)

\bibitem{Gomez2020RNC}
G{\'{o}}mez-Guti{\'{e}}rrez, D.: {On the design of nonautonomous fixed-time
  controllers with a predefined upper bound of the settling time}.
\newblock International Journal of Robust and Nonlinear Control
  \textbf{30}(10), 3871--3885 (2020).
\newblock \doi{10.1002/rnc.4976}.
\newblock \urlprefix\url{https://doi.org/10.1002/rnc.4976}

\bibitem{Jimenez2019}
Jimenez-Rodriguez, E., Munoz-Vazquez, A.J., Sanchez-Torres, J.D., Defoort, M.,
  Loukianov, A.G.: {A Lyapunov-Like Characterization of Predefined-Time
  Stability}.
\newblock IEEE Transactions on Automatic Control \textbf{65}(11), 4922--4927
  (2020).
\newblock \doi{10.1109/TAC.2020.2967555}

\bibitem{Khalil2002NonlinearSystems}
Khalil, H.K.: {Nonlinear systems}, third edn.
\newblock Prentice Hall (2002)

\bibitem{Liberzon2019CalculusTheory}
Liberzon, D.: {Calculus of Variations and Optimal Control Theory}.
\newblock Calculus of Variations and Optimal Control Theory  (2019)

\bibitem{Pal2020DesignTime}
Pal, A.K., Kamal, S., Nagar, S.K., Bandyopadhyay, B., Fridman, L.: {Design of
  controllers with arbitrary convergence time}.
\newblock Automatica \textbf{112} (2020).
\newblock \doi{10.1016/j.automatica.2019.108710}

\bibitem{Polyakov2012a}
Polyakov, A.: {Nonlinear feedback design for fixed-time stabilization of linear
  control systems}.
\newblock IEEE Transactions on Automatic Control \textbf{57}(8), 2106--2110
  (2012).
\newblock \doi{10.1109/TAC.2011.2179869}

\bibitem{Sanchez-Torres2018}
S{\'{a}}nchez-Torres, J.D., G{\'{o}}mez-Guti{\'{e}}rrez, D., L{\'{o}}pez, E.,
  Loukianov, A.G.: {A class of predefined-time stable dynamical systems}.
\newblock IMA Journal of Mathematical Control and Information \textbf{35}(1),
  I1--I29 (2018).
\newblock \doi{10.1093/imamci/dnx004}

\bibitem{Sanchez-Torres2020ASystems}
S{\'{a}}nchez-Torres, J.D., Mu{\~{n}}oz-V{\'{a}}zquez, A.J., Defoort, M.,
  Jim{\'{e}}nez-Rodr{\'{i}}guez, E., Loukianov, A.G.: {A class of
  predefined-time controllers for uncertain second-order systems}.
\newblock European Journal of Control \textbf{53}, 52--58 (2020).
\newblock \doi{10.1016/j.ejcon.2019.10.003}

\bibitem{Seeber2020ConvergenceControllers}
Seeber, R.: {Convergence Time Bounds for a Family of Second-Order Homogeneous
  State-Feedback Controllers}.
\newblock IEEE Control Systems Letters \textbf{4}(4), 1018--1023 (2020).
\newblock \doi{10.1109/LCSYS.2020.2998673}

\bibitem{ShihongDing2015NewControllers}
{Shihong Ding}, Levant, A., Li, S.: {New families of high-order sliding-mode
  controllers}.
\newblock In: 2015 54th IEEE Conference on Decision and Control (CDC), pp.
  4752--4757. 2015 54th IEEE Conference on Decision and Control (CDC) (2015)

\bibitem{Shtessel2014ObservationObservers}
Shtessel, Y., Edwards, C., Fridman, L., Levant, A.: {Observation and
  Identification via HOSM Observers}.
\newblock In: Sliding Mode Control and Observation, Control Engineering, pp.
  251--290. Birkh{\"{a}}user, New York, NY, New York, NY (2014).
\newblock \doi{10.1007/978-0-8176-4893-0{\_}7}

\bibitem{Song2017Time-varyingTime}
Song, Y., Wang, Y., Holloway, J., Krstic, M.: {Time-varying feedback for
  regulation of normal-form nonlinear systems in prescribed finite time}.
\newblock Automatica \textbf{83}, 243--251 (2017).
\newblock \doi{10.1016/J.AUTOMATICA.2017.06.008}

\bibitem{Song2019Time-varyingTime}
Song, Y., Wang, Y., Krstic, M.: {Time-varying feedback for stabilization in
  prescribed finite time}.
\newblock International Journal of Robust and Nonlinear Control \textbf{29}(3),
  618--633 (2019).
\newblock \doi{10.1002/rnc.4084}

\bibitem{Sontag2008InputResults}
Sontag, E.D.: {Input to State Stability: Basic Concepts and Results}.
\newblock Lecture Notes in Mathematics \textbf{1932}, 163--220 (2008)

\bibitem{Tabatabaeipour2014CalculationAnalysis}
Tabatabaeipour, S.M., Blanke, M.: {Calculation of Critical Fault Recovery Time
  for Nonlinear Systems based on Region of Attraction Analysis}.
\newblock IFAC Proceedings Volumes \textbf{47}(3), 6741--6746 (2014).
\newblock \doi{10.3182/20140824-6-ZA-1003.01418}.
\newblock
  \urlprefix\url{https://linkinghub.elsevier.com/retrieve/pii/S1474667016426716}

\bibitem{Utkin1992}
Utkin, V.I.: {Sliding Modes in Control and Optimization}.
\newblock Springer Verlag (1992).
\newblock \doi{10.1007/978-3-642-84379-2}

\bibitem{Zarchan2012}
Zarchan, P.: {Tactical and strategic missile guidance}.
\newblock American Institute of Aeronautics and Astronautics, Inc. (2012)

\bibitem{Zimenko2018}
Zimenko, K., Polyakov, A., Efimov, D., Perruquetti, W.: {On simple scheme of
  finite/fixed-time control design}.
\newblock International Journal of Control \textbf{93}(6), 1353--1361 (2020).
\newblock \doi{10.1080/00207179.2018.1506889}

\end{thebibliography}

\end{document}